\documentclass[12pt]{article}
\usepackage[margin=1in]{geometry}

\usepackage{mathptmx}

\usepackage{amsthm,amsmath,amsfonts,amssymb}
\usepackage[authoryear,round]{natbib}
\usepackage{microtype}
\usepackage{graphicx}
\providecommand{\alttext}[1]{}
\usepackage{float}
\usepackage{placeins}
\usepackage{needspace}
\usepackage{enumitem}
\usepackage{mathtools}
\usepackage{booktabs}
\usepackage{longtable}
\usepackage{xcolor}
\usepackage{tikz}
\usetikzlibrary{arrows.meta,calc,positioning,decorations.pathreplacing}
\usepackage[colorlinks,citecolor=blue,urlcolor=blue,filecolor=blue,linkcolor=blue,
         hypertexnames=false,plainpages=false]{hyperref}

\setlist*[enumerate]{label=(\roman*)}

\theoremstyle{plain}
\newtheorem{theorem}{Theorem}
\newtheorem{proposition}{Proposition}
\newtheorem{corollary}{Corollary}
\newtheorem{lemma}{Lemma}
\theoremstyle{definition}
\newtheorem{definition}{Definition}
\newtheorem{assumption}{Assumption}
\newtheorem{remark}{Remark}

\renewcommand{\qed}{\hfill\ensuremath{\square}}

\newdimen\figwidth

\renewcommand{\paragraph}[1]{\par\medskip\noindent\textit{#1}\enspace}

\DeclareMathOperator*{\argmin}{arg\,min}
\DeclareMathOperator*{\argmax}{arg\,max}
\DeclareMathOperator{\TV}{TV}
\DeclareMathOperator{\tr}{tr}
\DeclareMathOperator{\diag}{diag}

\newcommand{\Pcal}{\mathcal{P}}
\newcommand{\Ccal}{\mathcal{C}}
\newcommand{\Xcal}{\mathcal{X}}
\newcommand{\Acal}{\mathcal{A}}

\newcommand{\Ecal}{\mathcal{E}}

\newcommand{\Nor}{\mathcal{N}}
\newcommand{\R}{\mathbb{R}}
\newcommand{\Ebb}{\mathbb{E}}
\newcommand{\Pbb}{\mathbb{P}}
\newcommand{\ind}{\mathbf{1}}
\newcommand{\KL}{D_{\mathrm{KL}}}
\newcommand{\whatP}{\widehat{P}_n}

\newcommand{\DeltaK}{\Delta^{k_r-1}}

\title{Bayesian Prediction under Moment Conditioning}

\author{Nicholas G. Polson\thanks{Booth School of Business, University of
Chicago, Chicago, Illinois 60637, U.S.A. \texttt{ngp@chicagobooth.edu}.}
\and Daniel Zantedeschi\thanks{Muma College of Business, University of South
Florida, Tampa, Florida 33620, U.S.A. \texttt{danielz@usf.edu}.}}

\date{\today}

\begin{document}

\maketitle

\begin{abstract}
Moment restrictions specify a class of laws, not a predictive model.  We
obtain one by conditioning an independent sample from a reference law on
its empirical moments, and define prediction as the law of a fixed block
selected from that conditioned ensemble.  On a finite partition this law is
an exact mixture over empirical types.  Under exact feasibility and lattice
regularity, the mixing law has a Gaussian limit on the feasible tangent
space, governed by the reduced Hessian, and the selected block approaches
independent sampling from the Kullback--Leibler projection.  A separate
finite-sample bound gives the same product limit for general real-valued
restrictions without lattice assumptions.  Refinement recovers the
projection on the original sample space.  Parameterizing the projected
family produces a predictive product criterion with a local
inverse-covariance expansion, connecting the construction to generalized
method of moments.
\end{abstract}


\thispagestyle{empty}
\newpage
\setcounter{page}{1}

\section{Introduction}\label{sec:intro}

Many statistical models restrict moments without specifying a complete
likelihood.  Survey calibration supplies population moments, portfolio
problems impose restrictions on return and risk, and generalized estimating
equations determine features of a law without determining the law itself.
Because such restrictions describe a class of admissible distributions,
prediction also requires a reference law, denoted by $Q$.
Here the Bayesian operation is conditioning at the predictive-law level:
$Q$ is fixed, and the empirical moment event updates its product law; no
prior over an unknown sampling distribution is introduced.

Among the laws dominated by $Q$ and satisfying the moments, the
Kullback--Leibler minimizer $P^*$ has the familiar Gibbs form
\citep{jaynes1957information,csiszar1975}.  This variational principle
selects the admissible law closest to $Q$; by itself, however, it gives no
probabilistic reason for observations to follow that law.  We ask whether
one emerges from conditioning: after an independent ensemble drawn from
$Q$ is conditioned on its empirical moments, does a fixed selected block
behave as independent observations from $P^*$, and at what rate?  This is
exchangeable-block prediction within the conditioned ensemble, rather than
a sequential forecast based on an observed past.

We first answer this question on a finite measurable partition.  The moment
restrictions then form an affine slice of the probability simplex.
Conditional on the moment event, the law of a selected block is an exact
mixture of hypergeometric laws over the feasible empirical types.
Approximation reduces accordingly to two steps: locating the types that
carry conditional mass, and comparing sampling without replacement from
each such type with independent sampling from its cell probabilities.

The sharper result requires an arithmetic compatibility condition.  When
the moment equations admit exactly feasible empirical cell frequencies and
a local lattice regularity condition holds, the type weights have a Gaussian
limit on the feasible tangent space (Theorem~\ref{thm:gdf}).  Combining this
limit with the finite-population correction gives the block-collapse bound
of Theorem~\ref{thm:collapse}.  General real-valued moments need not admit
such frequencies.  Conditioning instead on a shrinking moment window gives
a less sharp finite-sample bound, but requires no lattice assumption
(Theorem~\ref{thm:window}).

These statements are initially finite-partition results.  A fixed-sample
Le~Cam comparison controls the loss from observing the partition rather
than the original sample, and refinement recovers the information
projection on the original sample space (\S\,\ref{sec:lecam}).
Parametrizing the projected family then produces a conditional likelihood
with a local generalized method-of-moments form
(\S\,\ref{sec:betel}).  The two operations
involve different notions of information: the Hessian restricted to the
constraint tangent space governs fluctuations of cell probabilities,
whereas parameter information appears only after differentiating the
parametrized family.

Thus the argument has two branches.  The predictive branch passes from the
conditioned reference ensemble to its exact type mixture, then to the
product law $(p_r^*)^{\otimes m}$ on a chart and finally, by lifting and
refinement, to $(P^*)^{\otimes m}$ on the original sample space.  The
inferential branch parameterizes the projected laws
$\theta\mapsto p_{\theta,r}^*$ and uses their product scores, whose local
expansion yields the inverse cell-moment covariance.

The probabilistic argument draws on exchangeability
\citep{deFinetti37,DiaconisFreedman1980,Aldous1985,DiaconisFreedman1987},
large deviations \citep{sanov1957,dembo1998large,Ellis2006,Lanford1973},
and Gibbs conditioning \citep{vancampenhout1981,csiszar1984}.  Its
inferential consequences connect with empirical likelihood
\citep{owen1988empirical,owen2001bel}, generalized method of moments
\citep{hansen1982gmm}, exponential tilting
\citep{qinlawless1994,schennach2007el,chibsimoni2018,chibsimoni2022,
Lazar2021}, and information geometry
\citep{AmariNagaoka2000,KassVos1997}.  The closest construction is that of
\citet{BornnShephardSolgi2019}, who take the constrained manifold as
primitive; here the finite-dimensional constraint set is derived from the
empirical measure through a partition.

The companion paper \citep{polson2025definettisanov} instead places a prior
over reference laws and studies asymptotic prediction after averaging over
them.  The present paper fixes $Q$ in order to retain the finite-sample
mixture representation, obtain a local limit theorem and quantify predictive
approximation on a partition.  The two approaches therefore address
complementary aspects of moment-conditioned prediction.

Section~\ref{sec:setup} formulates the constrained problem, and
\S\S\,\ref{sec:main}--\ref{sec:lecam} establish predictive concentration
on a partition and connect it to the original experiment.
\S\S\,\ref{sec:connections}--\ref{sec:misspec} develop the asymptotic,
parametric and misspecification consequences; \S\,\ref{sec:simulation}
gives a computational illustration, and \S\,\ref{sec:discussion}
concludes.  Extended proofs appear in the Appendix, which opens with a table of notation
(p.~\pageref{sec:notation}).

\section{The Constrained Predictive Problem}\label{sec:setup}

\subsection{Ambient Measurable-Space Setup}\label{sec:ambient}

Let $(\Xcal,\Acal)$ be a measurable space and let $Q\in\Pcal(\Xcal)$
be a baseline probability measure.  Let $h:\Xcal\to\R^d$ be a
measurable moment map with $h\in L^1(Q;\R^d)$, so that the conditional
cell moments used below are finite.  For a target value $\alpha\in\R^d$, define the
admissible class
\begin{equation}\label{eq:admissible}
\Ccal(\alpha)=\left\{P\in\Pcal(\Xcal):P\ll Q,\;\int h\,dP=\alpha\right\}.
\end{equation}
The absolute continuity condition $P\ll Q$ ensures that the
constrained problem has a well-defined density $dP/dQ$ and that the
Csisz\'ar--Bregman projection theorem applies.

\begin{assumption}\label{ass:regularity}
The log-moment-generating function
$\psi(\lambda)=\log\int e^{\lambda^\top h(x)}\,Q(dx)$ is essentially
smooth and steep, with $\Theta_\psi=\mathrm{int}\,\operatorname{dom}\psi\neq
\varnothing$, and the target satisfies
$\alpha\in\nabla\psi(\Theta_\psi)$, so that $\nabla\psi(\lambda^*)=\alpha$
has a solution $\lambda^*\in\Theta_\psi$; moreover
$\nabla^2\psi(\lambda^*)$ is positive definite.
\end{assumption}

The condition $\alpha\in\nabla\psi(\Theta_\psi)$ is the analogue of the
empirical-likelihood convex-hull condition
\citep{owen2001bel,schennach2005}; on a finite chart it becomes the
relative-interior condition of Definition~\ref{def:adequate}, which
yields a strictly positive projection (Lemma~A1).

\begin{theorem}[Measure-level projection]\label{thm:projection}
Under Assumption~\ref{ass:regularity}, the information projection
\[
P^*=\argmin_{P\in\Ccal(\alpha)}\KL(P\|Q)
\]
exists, is unique, and has Gibbs form
\begin{equation}\label{eq:gibbs}
\frac{dP^*}{dQ}(x)=\exp\{\lambda^{*\top}h(x)-\psi(\lambda^*)\},
\end{equation}
where $\lambda^*\in\R^d$ is the unique solution of
$\nabla\psi(\lambda^*)=\alpha$.  The value of the divergence is
$\KL(P^*\|Q)=\lambda^{*\top}\alpha-\psi(\lambda^*)$.
\end{theorem}

The theorem requires no finiteness assumption on $\Xcal$.  It is the
classical I-projection theorem \citep{csiszar1975}; the Pythagorean identity
$\KL(P\|Q)-\KL(P^*\|Q)=\KL(P\|P^*)$ gives the key optimality step.
Moreover, $\nabla^2\psi(\lambda^*)=\mathrm{Cov}_{P^*}\{h(X)\}$.
The proof and functional-analytic details are in the Appendix, \S\,A1.  Our analysis uses the result on the finite empirical
coordinates introduced in \S\,\ref{sec:partition}--\ref{sec:geometry}.

\begin{remark}[Choice of baseline measure]\label{rem:baseline}
The choice of $Q$ is part of the model: it records what is known before the
moment restriction is imposed.  It might be a reference member of a
parametric family, a kernel estimate obtained from independent training data
and then held fixed, or a maximum-entropy law.  Taking $Q=\whatP$ from the same
sample is different.  It produces a data-dependent pseudo-likelihood chart
and is not covered by the fixed-$Q$ theorem
(\S\,\ref{sec:betel}).  Since both $P^*$ and its curvature inherit this
choice, \S\,\ref{sec:misspec} treats misspecification of the baseline
explicitly.
\end{remark}

\subsection{Permutation Invariance and Empirical Reduction}\label{sec:permutation}

For $\Xcal$-valued $X_1,\ldots,X_n$, neither the empirical moment
$M_n=n^{-1}\sum_i h(X_i)$ nor the event that constrains it depends on the
order of the observations.  The analysis below uses a finite partition, on
which every permutation orbit is represented by the empirical type.
Consequently the conditional block law can be decomposed by type, without
choosing an ordering of the sample.  The ambient empirical measure
$\whatP=n^{-1}\sum_i\delta_{X_i}$ is the corresponding order-free
representation before discretization.

\begin{remark}[Exchangeability versus permutation invariance]
Two symmetries are easily confused here.  Permutation invariance of the
empirical constraint is algebraic, whereas exchangeability of the
generating sequence is probabilistic.  Exchangeability justifies selecting
an arbitrary block; the geometry itself comes from the empirical measure
and its constraint.
\end{remark}

\subsection{Finite-Partition Realization}\label{sec:partition}

Although $\Xcal$ may be infinite, an empirical measure based on $n$
observations has at most $n$ support points.  A finite measurable partition
$\Pi_r=\{A_1^{(r)},\ldots,A_{k_r}^{(r)}\}$ offers coordinates more stable
than the moving support of the empirical measure: a law is represented by
its cell probabilities in
$\Delta^{k_r-1}$.  At a fixed resolution $r$, the constrained problem is
therefore finite-dimensional and its tangent spaces and Hessians can be
written down directly.

\begin{definition}[Moment-adequate partition]\label{def:adequate}
A partition $\Pi_r$ is \emph{moment-adequate} for $(h,\alpha)$ if:
\begin{enumerate}
\item $q_{r,j}=Q(A_j^{(r)})>0$ for all $j=1,\ldots,k_r$;
\item the augmented constraint matrix
  $\widetilde A_r=(\ind\;A_r^\top)^\top$ has
  full row rank $d+1$, where $A_r=[h_r(1)\;\cdots\;h_r(k_r)]$ and
  $h_r(j)=\Ebb_Q\{h(X)\mid X\in A_j^{(r)}\}$;
\item the target lies in the relative interior,
  $\alpha\in\mathrm{ri}\,\mathrm{conv}\{h_r(1),\ldots,h_r(k_r)\}$.
\end{enumerate}
These are properties of a single partition; condition~(iii) guarantees a
strictly positive finite I-projection.  The refinement behaviour of a
partition \emph{sequence} is imposed separately in \S\,\ref{sec:lecam}
(Assumption~\ref{ass:refine}).
\end{definition}

Fix a moment-adequate partition.  Any probability measure $P$ then induces
the discrete law $p^{(r)}(j)=P(A_j^{(r)})$ on
$\{1,\ldots,k_r\}$, while the baseline induces
$q_{r,j}=Q(A_j^{(r)})$.  In these coordinates the constraint set is
\[
\Ccal_r(\alpha)=\left\{p\in\Delta^{k_r-1}:
\sum_{j=1}^{k_r}p_j\,h_r(j)=\alpha\right\}.
\]
This is the natural discretization because lifting the cell law $p$ with
the baseline conditionals $Q(\,\cdot\mid A_j^{(r)})$ gives an ambient law
whose moment is exactly $A_rp$; the approximation to the original moment
map is controlled under refinement in \S\,\ref{sec:lecam}.
Its information projection,
$p_r^*=\argmin_{p\in\Ccal_r(\alpha)}\KL(p\|q_r)$, again has Gibbs form:
\begin{equation}\label{eq:discrete-gibbs}
p_{r,j}^*=q_{r,j}\exp\{\lambda_r^{*\top}h_r(j)-\psi_r(\lambda_r^*)\},
\quad
\psi_r(\lambda)=\log\sum_{j=1}^{k_r}q_{r,j}\,e^{\lambda^\top h_r(j)}.
\end{equation}

Here $P$ denotes a measure on $(\Xcal,\Acal)$ and $p$ its partition-induced
vector, $p_j=P(A_j^{(r)})$.  Throughout, $\|\cdot\|_2$ and
$\|\cdot\|_1$ are the Euclidean and $\ell_1$ norms on $\R^{k_r}$, while
unsubscripted $\|\cdot\|$ is the Euclidean norm on $\R^d$.

\subsection{Constraint Geometry and Curvature}\label{sec:geometry}

The tangent space at $p_r^*$ contains exactly those perturbations that
preserve total mass and the prescribed moments:
\begin{equation}\label{eq:tangent}
T_r^*=\left\{v\in\R^{k_r}:\ind^\top v=0,\;
\sum_{j=1}^{k_r}v_j\,h_r(j)=0\right\},
\end{equation}
generically with dimension $k_r-1-d$.  With cell moment matrix
$A_r=[h_r(1)\;\cdots\;h_r(k_r)]$ and augmented constraint matrix
$\widetilde A_r=(\ind\;A_r^\top)^\top$ of full row rank, the orthogonal
projector onto $T_r^*$ is
$\Pi_{T_r^*}=I_{k_r}-\widetilde A_r^\top
(\widetilde A_r\widetilde A_r^\top)^{-1}\widetilde A_r$.  If
$V_r\in\R^{k_r\times(k_r-1-d)}$ is a matrix whose columns form an
orthonormal basis for $T_r^*$, so that $\Pi_{T_r^*}=V_rV_r^\top$,
the \emph{reduced Hessian} is
\begin{equation}\label{eq:hessian}
H_r^*=V_r^\top\diag(1/p_r^*)\,V_r.
\end{equation}
Here $\diag(1/p_r^*)$ denotes the diagonal matrix with $j$th entry
$1/p_{r,j}^*$; thus $H_r^*$ is the Kullback--Leibler Hessian, evaluated
at the projection and restricted to feasible directions.  It is a concrete
$(k_r-1-d)\times(k_r-1-d)$ matrix obtained from $p_r^*$ and a choice of
tangent basis.

\begin{proposition}[Quadratic expansion]\label{prop:quadratic}
For any $p\in\Ccal_r(\alpha)\cap(\Delta^{k_r-1})^\circ$ with
$p=p_r^*+V_r u$ for $u\in\R^{k_r-1-d}$ and
$\|u\|_2\le p_{\min,r}^*/2$,
\[
\KL(p\|q_r)=\KL(p_r^*\|q_r)+\tfrac{1}{2}u^\top H_r^* u+R_3(u),
\]
where $|R_3(u)|\le 2\|u\|_2^3/\{3\,(p_{\min,r}^*)^2\}$ and
$p_{\min,r}^*=\min_{j:q_{r,j}>0}p_{r,j}^*$; the radius condition
keeps every coordinate of $p$ above $p_{\min,r}^*/2$ along the
segment from $p_r^*$ to $p$.
\end{proposition}

The linear term vanishes because the Gibbs gradient is orthogonal to
$T_r^*$; restricting the Hessian $\diag(1/p_r^*)$ gives $H_r^*$.
The third-derivative calculation and remainder bound are given in the
Appendix, \S\,A2.

The role of $H_r^*$ is local: its spectrum governs concentration of
probability vectors on this chart, and hence predictive collapse, while
the parameter rates of \S\,\ref{sec:parametric} involve different
information matrices.  For the chart itself, the useful
quantities are the minimal curvature
$\lambda_{\min}(H_r^*)>0$, the upper bound
$\lambda_{\max}(H_r^*)\le 1/p_{\min,r}^*$, and the condition number
$\kappa_r=\lambda_{\max}(H_r^*)/\lambda_{\min}(H_r^*)$, which measures
local isotropy (Fig.~\ref{fig:geometry}).

\begin{figure}[htbp]
\centering
\begin{tikzpicture}[scale=2.0]
\coordinate (v1) at (0,0);
\coordinate (v2) at (2,0);
\coordinate (v3) at (1,1.732);
\draw[thick] (v1) -- (v2) -- (v3) -- cycle;
\node[below left,font=\small] at (v1) {$e_1$};
\node[below right,font=\small] at (v2) {$e_{k_r}$};
\node[above,font=\small] at (v3) {$e_2$};
\node[font=\footnotesize,gray] at (1.65,0.15) {$\DeltaK$};
\draw[very thick,black] (0.28,0.795) -- (1.70,0.706);
\node[black,font=\footnotesize,above right] at (0.27,0.80)
  {$\Ccal_r(\alpha)$};
\draw[black!40,thin,dotted] (0.95,0.5) ellipse (0.25 and 0.14);
\draw[black!40,thin,dotted] (0.95,0.5) ellipse (0.45 and 0.252);
\draw[-{Stealth[length=3.5pt]},thick,black]
  (0.95,0.5) -- (1.0,0.75);
\fill[black!55] (0.95,0.5) circle (1.2pt);
\node[black!55,font=\footnotesize,below] at (0.95,0.47) {$q_r$};
\fill[black] (1.0,0.75) circle (1.4pt);
\node[black,font=\footnotesize,above right] at (1.03,0.76) {$p_r^*$};
\end{tikzpicture}
\caption{Constraint geometry on $\DeltaK$.  The feasible set
$\Ccal_r(\alpha)$ is an affine slice of the simplex, parallel to the tangent
space $T_r^*$.  The baseline $q_r$ projects by Kullback--Leibler minimization to $p_r^*$;
the first touching divergence contour is tangent to the slice there, and the
reduced Hessian $H_r^*$ records curvature along it.}\label{fig:geometry}
\alttext{Simplex diagram with a straight affine constraint slice.  An arrow
projects the baseline to the touching point of an elliptical divergence
contour on the slice.}
\end{figure}

\subsection{Computation of Curvature}\label{sec:computation}

Compute $\lambda_r^*$ from $\nabla\psi_r(\lambda_r^*)=\alpha$, hence
$p_r^*$ by~\eqref{eq:discrete-gibbs}; orthogonalize the null space of
$\widetilde A_r$ to obtain $V_r$; then form
$H_r^*=V_r^\top\diag(1/p_r^*)V_r$ and compute its eigenvalues.
By the Rayleigh quotient of $\diag(1/p_r^*)$ restricted to $T_r^*$,
all eigenvalues of $H_r^*$ lie between $1/p_{\max,r}^*$ and
$1/p_{\min,r}^*$; the numerical stability of
$\lambda_{\min}(H_r^*)$ under perturbation of
$p_r^*$ is discussed in the Appendix,
Lemma~A2.

\FloatBarrier
\section{Predictive Concentration}\label{sec:main}

\subsection{Predictive Representation}\label{sec:predictive}

We define prediction by selecting a block from the moment-conditioned
ensemble itself.  Write
$Z_i=T_r(X_i)\in\{1,\ldots,k_r\}$ for the cell labels and take, as a
finite-chart baseline, $Z_1,\ldots,Z_n$ to be independent draws from
$q_r$.  Let $\whatP^{(r)}=n^{-1}\sum_{i=1}^n e_{Z_i}$ be
the empirical type and
\begin{equation}\label{eq:feasible}
E_{n,r}=\bigl\{p\in\Pcal_n(\{1,\ldots,k_r\}):
\|A_r p-\alpha\|\le\eta_n\bigr\}
\end{equation}
the set of approximately feasible types.  Here
$\Pcal_n(\{1,\ldots,k_r\})=\{p\in\Delta^{k_r-1}:np_j\in\mathbb{Z}_{\ge0}
\text{ for all }j\}$ is the set of $n$-types, $A_r$ is the cell-moment
matrix of \S\,\ref{sec:geometry}, and the tolerance $\eta_n\ge0$ determines
how closely a realizable type must satisfy the target moment.

The sharp result must accommodate a small arithmetic obstruction: at a
given sample size, the continuous constraint surface may pass between all
empirical types.  We say that \emph{lattice regularity} holds when the surface
$\{p:A_r p=\alpha,\ \ind^\top p=1\}$ is rational and, along a cofinal set
$\mathcal{N}_r$, its exactly feasible types form a translate of a
rank-$s_r$ lattice, $s_r=k_r-1-d$, with mesh $O(n^{-1})$.  We also require
one such type within $O_r(n^{-1})$ of $p_r^*$.  Rational $h_r(j)$ and
$\alpha$, together with a full-rank integer kernel, are sufficient.

The condition matters because Gibbs coordinates are usually irrational,
so there is no reason for $p_r^*$ itself to be an $n$-type.  The sample
can instead realize the discrete minimizer
$\bar p_{n,r}=\argmin\{\KL(p\|q_r):p\in\Pcal_n,\ A_r p=\alpha\}$.
Under lattice regularity this lies within $O_r(n^{-1})$ of $p_r^*$, and
because every feasible type belongs to $p_r^*+T_r^*$, the parametrization
$\widehat p=p_r^*+V_r u/\sqrt n$ is exact and the Gaussian is centred at
$p_r^*$, whether or not that point belongs to the lattice.  The more
forgiving case $\eta_n>0$ is treated in Theorem~\ref{thm:window} and
Remark~\ref{rem:thickening}.

For a block size $1\le m\le n$, the constrained predictive law is the
conditional law of the first $m$ labels given the moment event,
\begin{equation}\label{eq:predictive}
\mu_{n,m}^{(r)}(A)
=\Pbb\bigl\{(Z_1,\ldots,Z_m)\in A \bigm| \whatP^{(r)}\in E_{n,r}\bigr\},
\qquad A\subseteq\{1,\ldots,k_r\}^m .
\end{equation}
The phrase ``first $m$'' is only notation: because the ensemble is
exchangeable and the conditioning event depends only on its type, these
coordinates may be replaced by any fixed $m$ coordinates, or by an
$m$-subset chosen before the labels are seen.  Thus
$\mu_{n,m}^{(r)}$ describes a \emph{selected exchangeable block} from the
conditioned ensemble, rather than a forecast formed from an observed past
and an independent future.

The value of the definition appears on conditioning once more, now on the
empirical type.  For a type $p$ with counts $np_j$, let
\begin{equation}\label{eq:typeweight}
w_{n,r}(p)\ \propto\ \ind\{p\in E_{n,r}\}\,
\dbinom{n}{np_1,\ldots,np_{k_r}}\prod_{j}q_{r,j}^{\,np_j},
\qquad \textstyle\sum_{p}w_{n,r}(p)=1,
\end{equation}
and, for a block $z_{1:m}$ with $c_j(z_{1:m})$ occurrences of cell $j$,
the within-type sampling-without-replacement law
\begin{equation}\label{eq:hypergeom}
\mathrm{HG}_{n,p,m}(z_{1:m})
=\frac{\prod_{j=1}^{k_r}(np_j)_{c_j(z_{1:m})}}{(n)_m},
\qquad (a)_b=a(a-1)\cdots(a-b+1).
\end{equation}
Given the type, the selected coordinates are draws without replacement;
averaging this law over the possible types therefore gives the
\emph{exact} identity
\begin{equation}\label{eq:mixture}
\mu_{n,m}^{(r)}=\sum_{p\in E_{n,r}}
\mathrm{HG}_{n,p,m}\,w_{n,r}(p).
\end{equation}
The weights are multinomial because the finite-chart baseline draws are
independent and identically distributed;
a different exchangeable baseline would contribute its own mixing law
over types.  The localization and collapse results below use the
multinomial weights and therefore retain the independent baseline
assumption.  Although \citet{DiaconisFreedman1980} provide the surrounding
finite-exchangeability theory, \eqref{eq:mixture} appeals to no
approximation theorem: it is an exact consequence of conditioning.

Approximation enters at the next step, where the collision coupling shows
that, within a fixed type, the block law is close to $p^{\otimes m}$:
\begin{equation}\label{eq:collision}
\|\mathrm{HG}_{n,p,m}-p^{\otimes m}\|_{\TV}
\le 1-\frac{(n)_m}{n^m}\le\frac{m(m-1)}{2n},
\end{equation}
uniformly in $p$.  Indeed, sampling with and without replacement may be
coupled through the same draws until an index is repeated.  It then remains
only to locate the random type; \S\,\ref{sec:localization} shows that
its weight concentrates near $p_r^*$, and the two approximations together
yield the product limit $(p_r^*)^{\otimes m}$.

\subsection{Localization}\label{sec:localization}

The concentration argument separates a shrinking neighbourhood of
$p_r^*$, where the Kullback--Leibler divergence admits a quadratic
expansion, from its
complement, where a large-deviation bound suffices.  This good--bad
decomposition is in the spirit of Lanford's analysis
\citep{Lanford1973}.

Define the curvature radius
$\rho_n=\{4\log n/\lambda_{\min}(H_r^*)\}^{1/2}$ and partition the
feasible types into
\begin{equation}\label{eq:localization}
\mathcal{G}_{n,r}
=\{\widehat{p}\in E_{n,r}:\|\widehat{p}-p_r^*\|_2\le\rho_n/\sqrt{n}\},
\quad
\mathcal{B}_{n,r}=E_{n,r}\setminus\mathcal{G}_{n,r}.
\end{equation}
The choice of $\rho_n$ is calibrated to the quadratic Kullback--Leibler
term: crossing the
boundary of $\mathcal{G}_{n,r}$ costs at least $\log n/n$ in divergence.
After multiplication by $n$ in the multinomial exponent, this produces a
factor of at most $n^{-1}$.

\begin{lemma}[Exponential suppression of bad types]\label{lem:tail}
Let the exact-feasibility and lattice-regularity conditions of
\S\,\ref{sec:predictive} hold, along $n\in\mathcal{N}_r$.  With $\rho_n$
as above,
\[
\sum_{\widehat{p}\in\mathcal{B}_{n,r}}
\Pbb(\whatP^{(r)}=\widehat{p}\mid \whatP^{(r)}\in E_{n,r})
=O\{n^{-1}(\log n)^{s_r/2}\}.
\]
\end{lemma}

The proof splits the bad set into a local annulus, controlled by the
quadratic expansion, and a remote region with a fixed divergence gap.  A
lattice-sum comparison then normalizes the local tail; see the
Appendix, \S\,A3.

\begin{lemma}[Gaussian control of good types]\label{lem:gaussian}
Under the conditions of Lemma~\ref{lem:tail}, within
$\mathcal{G}_{n,r}$ the conditional weight of each type is
\[
\Pbb(\whatP^{(r)}=\widehat{p}\mid \whatP^{(r)}\in E_{n,r})
\propto\exp\left\{-\tfrac{1}{2}u^\top H_r^*\,u\right\}
\cdot\{1+O(\rho_n^3 n^{-1/2}+n^{-1})\},
\]
where $\widehat{p}=p_r^*+V_r u/\sqrt{n}$ with
$u\in\R^{s_r}$, $V_ru\in T_r^*$, and
$\|u\|_2\le\rho_n$.
\end{lemma}

Stirling's formula writes the multinomial type probability as a slowly
varying prefactor times $\exp\{-n\KL(\widehat p\|q_r)\}$.  Substituting the
quadratic expansion at $p_r^*$ yields the displayed Gaussian weight; the
uniform prefactor and lattice errors are controlled in the Appendix, \S\,A3.

\subsection{Main Results}\label{sec:results}

The next three results separate distributional approximation from
predictive concentration.  Under exact lattice feasibility,
Theorem~\ref{thm:gdf} approximates the conditional law by a Gaussian mixture
along the feasible tangent directions, and Theorem~\ref{thm:collapse}
integrates out those fluctuations to bound the distance to
$(p_r^*)^{\otimes m}$.  Theorem~\ref{thm:window} retains product collapse
for generic real-valued moments, using a shrinking window instead of a
lattice local limit.

\begin{theorem}[Gaussian localization on the constraint chart]\label{thm:gdf}
Fix a moment-adequate partition $\Pi_r$, let
Assumption~\ref{ass:regularity} and the lattice regularity of
\S\,\ref{sec:predictive} hold, and work in the exact-feasibility case
$\eta_n=0$ along the cofinal set $n\in\mathcal{N}_r$.  For $1\le m\le n$
the constrained predictive law satisfies
\begin{equation}\label{eq:gdf}
\mu_{n,m}^{(r)}(A)
=\int_{\|v\|_2\le\rho_n} \{p_r^*(v)\}^{\otimes m}(A)\;\varphi_{H_r^*}(v)\,dv
+R_{n,m,r}(A),
\end{equation}
where $p_r^*(v)=p_r^*+V_r v/\sqrt{n}$ and $\varphi_{H_r^*}$ is the
Gaussian density on $\R^{s_r}$ with precision $H_r^*$.  The remainder
collects four contributions, each with a constant depending only on the
fixed chart $r$: the local Taylor and Stirling error; the
lattice-to-integral (coarea) error; the finite-population error
$O(m^2/n)$ incurred when $\mathrm{HG}_{n,p,m}$ is replaced by
$p^{\otimes m}$ via~\eqref{eq:collision}; and an exponentially small
Gaussian-tail error.  For fixed $r$,
\[
\sup_A|R_{n,m,r}(A)|
\le C_r\bigl\{(1+m)\,n^{-1/2}(\log n)^{3/2}+m^2/n\bigr\}.
\]
\end{theorem}

\begin{proof}
Decompose the exact mixture~\eqref{eq:mixture} over $E_{n,r}$ into good
types $\mathcal{G}_{n,r}$ and bad types $\mathcal{B}_{n,r}$
(\S\,\ref{sec:localization}).  Lemma~\ref{lem:tail} bounds the bad-type
mass by $O_r\{n^{-1}(\log n)^{s_r/2}\}$, and
Lemma~\ref{lem:gaussian} replaces each good-type weight by a Gaussian
density with a fixed-$r$ relative error.  Within the good set the
collision bound~\eqref{eq:collision} replaces $\mathrm{HG}_{n,p,m}$ by
$p^{\otimes m}$ at cost $O(m^2/n)$; in the exact-feasibility case the
continuous constrained minimizer is $p_r^*$, and the discrete minimizer
lies within $O_r(n^{-1})$ of it, so no thickening term arises.  Passage
from the lattice sum to the
Gaussian integral is carried out in the affine simplex chart by a
coarea/covolume calculation whose covolume factor cancels in the
normalized weights, with a fixed-$r$ error; normalizing against the full
Gaussian integral on $\R^{s_r}$ adds an exponentially small tail.  Collecting the
four contributions gives $R_{n,m,r}$; the fixed-$r$ constants are worked
out in the Appendix, \S\S\,A3--A5.
\end{proof}

\begin{remark}[Gaussian measure on the constraint manifold]%
\label{rem:gaussian}
The representation replaces the conditional multinomial type law with a
Gaussian measure on the tangent space $T_r^*$: the constrained
predictive law is a Gaussian projected onto the
$(k_r\!-\!1\!-\!d)$-dimensional affine subspace
$\{p:\sum_j p_j h_r(j)=\alpha\}$, which is a reduced-rank Gaussian
with covariance determined by the reduced Hessian $H_r^*$, the
information matrix of the tangent-coordinate experiment on the chart;
\S\,\ref{sec:connections} separates this from parameter-space
information.
\end{remark}

\begin{theorem}[Exchangeable-block collapse]\label{thm:collapse}
Fix a moment-adequate partition $\Pi_r$.  Under
Assumption~\ref{ass:regularity} and the lattice regularity of
\S\,\ref{sec:predictive}, in the exact-feasibility case along
$n\in\mathcal{N}_r$ and for $1\le m\le n$,
\begin{equation}\label{eq:collapse}
\|\mu_{n,m}^{(r)}-(p_r^*)^{\otimes m}\|_{\TV}
\le C_r\left[
m\left(\frac{k_r\,\tr\{(H_r^*)^{-1}\}}{n}\right)^{1/2}
+\frac{m(m-1)}{n}
+\frac{(1+m)(\log n)^{3/2}}{\sqrt{n}}\right],
\end{equation}
with a constant $C_r$ depending only on the chart.  In particular, for
every fixed $r$ and every fixed block size $m$ the right-hand side tends
to $0$ as $n\to\infty$.
\end{theorem}

\begin{proof}
By Theorem~\ref{thm:gdf} and the triangle inequality the distance is at
most $\int_{\|v\|_2\le\rho_n}
\|p_r^*(v)^{\otimes m}-(p_r^*)^{\otimes m}\|_{\TV}\,
\varphi_{H_r^*}(v)\,dv+\sup_A|R_{n,m,r}(A)|$.  The product-coupling
inequality $\|P^{\otimes m}-Q^{\otimes m}\|_{\TV}\le m\|P-Q\|_{\TV}$ and
$\|p_r^*(v)-p_r^*\|_{\TV}=\|V_r v\|_1/(2\sqrt{n})
\le\sqrt{k_r}\,\|v\|_2/(2\sqrt{n})$, together with
$\int_{\|v\|\le\rho_n}\|v\|_2\varphi_{H_r^*}(v)\,dv
\le\Ebb_\varphi\{\|v\|_2\}\le\{\tr(H_r^{*-1})\}^{1/2}$ give the leading
term $m\{k_r\tr((H_r^*)^{-1})/n\}^{1/2}$, with no logarithmic factor.
The collision bound~\eqref{eq:collision} contributes $m(m-1)/n$, and the
Taylor, Stirling and coarea remainders of Theorem~\ref{thm:gdf} give the
final term.  Details are in the Appendix, \S\,A5.
\end{proof}

The leading term is the Gaussian tangential spread, with
$k_r^{1/2}$ converting Euclidean fluctuation to total variation; on a
balanced chart it is $O\{m(k_r/n)^{1/2}\}$.  The second term is the
finite-population collision correction, while the last collects the
localization, Taylor and lattice errors.  Here
$\rho_n^2=4\log n/\lambda_{\min}(H_r^*)$ makes the bad-type tail
$O_r\{n^{-1}(\log n)^{s_r/2}\}$.

Predictive collapse is a finite-chart, quantitative form of the Gibbs
conditioning principle \citep{vancampenhout1981,csiszar1984}:
conditionally on the moment event, a selected exchangeable block becomes
asymptotically independent with common law the information projection;
Theorem~\ref{thm:collapse} gives a finite-$n$ bound on the chart, with
constants depending on the chart geometry.

The sharper bound rests on lattice regularity.  For generic real-valued
charts, where the constraint surface misses the type lattice, a coarser
bound holds without that hypothesis.

\begin{theorem}[Window-robust collapse]\label{thm:window}
Fix a moment-adequate partition $\Pi_r$, let
Assumption~\ref{ass:regularity} hold, and take a tolerance with
$\eta_n\to0$ and $n\eta_n\to\infty$; lattice regularity is \emph{not}
assumed.  Then for $1\le m\le n$,
\begin{equation}\label{eq:window}
\|\mu_{n,m}^{(r)}-(p_r^*)^{\otimes m}\|_{\TV}
\le C_r\Bigl(m\eta_n+m\sqrt{\tfrac{\log n}{n}}+\tfrac{m^2}{n}\Bigr).
\end{equation}
\end{theorem}

\begin{proof}
Let
$p_r^{\eta_n}=\argmin\{\KL(p\|q_r):\|A_rp-\alpha\|\le\eta_n\}$
be the continuous slab minimizer.  For any feasible $p$, first-order optimality
on the convex slab and the Bregman identity for $\KL$ give
\[
\KL(p\|q_r)-\KL(p_r^{\eta_n}\|q_r)\ \ge\ \KL(p\|p_r^{\eta_n})
\ \ge\ \tfrac12\|p-p_r^{\eta_n}\|_1^2 ,
\]
the last step by Pinsker's inequality.  With the polynomial bound on the
number of $n$-types this yields
$\Ebb\{\|\whatP^{(r)}-p_r^{\eta_n}\|_1\mid\whatP^{(r)}\in E_{n,r}\}
=O_r(\sqrt{\log n/n})$, and $\|p_r^{\eta_n}-p_r^*\|_1=O_r(\eta_n)$
(Appendix, \S\,A4).  Product coupling and the collision
bound~\eqref{eq:collision} give the three terms.  No tangent expansion and
no lattice structure are used, so the bound holds for generic charts, at
the cost of a logarithm.  Details are in the Appendix, \S\,A5.
\end{proof}

\begin{remark}[Exact slice versus generic window]\label{rem:thickening}
For generic real-valued cell moments, the shrinking window ensures that
the conditioning event contains realizable types.  Its continuous minimizer
satisfies $\|p_r^{\eta_n}-p_r^*\|_2=O_r(\eta_n)$, which produces the first
term in Theorem~\ref{thm:window}; choosing
$\sqrt n\,\eta_n\to0$ keeps that term below the stochastic scale.
Unlike exact feasibility, the window intersects several normal layers of
the type lattice.  A Gaussian description would therefore require a
two-scale local limit, tangent at scale $n^{-1/2}$ and normal at scale
$n^{-1}$; the Appendix, \S\,A4, gives the details and records
this extension as open.
\end{remark}

Theorem~\ref{thm:window} provides the predictive justification for the
induced likelihood of \S\,\ref{sec:betel} on ordinary real-valued
partitions, where the exact-lattice hypothesis of
Theorem~\ref{thm:collapse} does not apply.  It is the chart form, with
chart-dependent constants $C_r$, of the unpartitioned master inequality of
the companion paper \citep{polson2025definettisanov}.  These constants are
required in \S\,\ref{sec:betel} but are not supplied by the ambient bound.

\begin{remark}[Claim calibration]\label{rem:calibration}
Theorems~\ref{thm:gdf} and~\ref{thm:collapse} are fixed-chart
statements, with $r$-dependent constants.  The approximation to the
ambient construction enters only through refinement
(\S\,\ref{sec:lecam}): the projected chart law approaches the ambient
projected law at rate $O(m\omega_r)$.
\end{remark}

\section{The Le~Cam Bridge}\label{sec:lecam}

A Le~Cam deficiency argument at the experiment level, together with a
projected-law approximation, connects the finite-coordinate results of
\S\,\ref{sec:main} to the ambient measure-theoretic construction under
partition refinement.

\subsection{Experiment Comparison}

Let $\{P_\theta:\theta\in\Theta\}$ be a parametric family on
$(\Xcal,\Acal)$ with moment targets $\alpha(\theta)$, as used in
\S\,\ref{sec:parametric}, let
$\Ecal_n=\{P_\theta^{\otimes n}:\theta\in\Theta\}$ denote the
corresponding experiment on $\Xcal^n$, and let
$T_r:\Xcal\to\{1,\ldots,k_r\}$ be the partition map.  The discretized
experiment is
$\Ecal_{n,r}=\{(P_\theta\circ T_r^{-1})^{\otimes n}:\theta\in\Theta\}$.

Since $T_r$ is a deterministic function of the data, the original
experiment can always reproduce the discretized one: the partition map
$T_r^{\otimes n}$ serves as the Markov kernel, so
\begin{equation}\label{eq:easy}
\delta(\Ecal_n,\Ecal_{n,r})=0.
\end{equation}
Coarsening the original experiment requires no additional information.
The reverse deficiency $\delta(\Ecal_{n,r},\Ecal_n)$ measures how well the
discretized experiment can reconstruct the original.  A \emph{reverse
kernel}
$\kappa_{r,j}(\cdot)=Q(\cdot\mid A_j^{(r)})$, which maps each cell
$j$ to the baseline conditional distribution within it, produces the
lifted law
\[
\widetilde{P}_{\theta,r}(dx)
=\sum_{j=1}^{k_r}P_\theta(A_j^{(r)})\,\kappa_{r,j}(dx).
\]
\begin{assumption}[Refinement regularity]\label{ass:refine}
The partitions $\{\Pi_r\}_{r\ge1}$ form a refining sequence with
$\sigma(\bigcup_r\Pi_r)=\Acal$, and the cell oscillation
$\omega_r=\sup_j\sup_{x,x'\in A_j^{(r)}}\|h(x)-h(x')\|\to0$ as
$r\to\infty$; on an unbounded sample space this is imposed after
truncation, or replaced by a direct cell-density regularity condition
(Appendix, \S\,A7).
\end{assumption}

Under the hypotheses of Theorem~A1, including the direct uniform
within-cell regularity of the densities $f_\theta=dP_\theta/dQ$, the
Le~Cam distance satisfies
\begin{equation}\label{eq:lecam-rate}
\Delta(\Ecal_n,\Ecal_{n,r})
\le n\cdot\sup_{\theta\in\Theta_0}\TV(P_\theta,\widetilde{P}_{\theta,r})
=O(n\omega_r).
\end{equation}
In particular, $\Delta(\Ecal_n,\Ecal_{n,r})\to 0$ as $r\to\infty$ for
fixed $n$.  The factor $n$ makes this an auxiliary experiment-level bound,
proved as Theorem~A1 in the Appendix.  Corollary~\ref{cor:joint} gives
the required finite-to-ambient result by allowing the resolution $r=r_n$
to grow with $n$, so that discretization bias and finite-sample
concentration vanish together.

\subsection{Predictive Lifting under Refinement}

For prediction, the ambient and finite limits are connected through
their projected laws, not by a direct conditional-deficiency claim.  Let
$\widetilde P_r^*=\sum_j p_{r,j}^*\,\kappa_{r,j}$ be the lifted
projection.  A projection-approximation lemma (Appendix,
Lemma~A3) shows that, under bounded $\lambda^*$, nonsingular moment
covariance, Assumption~\ref{ass:refine}, and the local exponential-envelope
condition stated there,
$\TV(\widetilde P_r^*,P^*)=O(\omega_r)$, so the lifted block-predictive
limits $(P^*)^{\otimes m}$ and $(\widetilde P_r^*)^{\otimes m}$, both on
$\Xcal^m$, differ by $O(m\omega_r)$.
Together with Theorem~\ref{thm:window}, this establishes the
ambient-to-finite approximation through the projected laws
(Corollary~\ref{cor:joint}).

\begin{remark}[Heuristic partition calibration]\label{rem:partition-sizing}
Under \emph{uniform} balanced-chart assumptions, namely
$p_{\min,r}^*\asymp k_r^{-1}$ together with uniform control over $r$ of
$\lambda_{\min}(H_r^*)\asymp k_r$, of the cubic and lattice constants of
Theorems~\ref{thm:gdf}--\ref{thm:collapse}, and of the convergence
$\lambda_r^*\to\lambda^*$, the fixed-chart bounds can be balanced
formally.  Coarsening the partition degrades the projected-law
approximation by $O(m\omega_r)$; refining it inflates the concentration
term $m\{k_r\tr((H_r^*)^{-1})/n\}^{1/2}$ of Theorem~\ref{thm:collapse},
which on a balanced chart is $C_2\,m(k_r/n)^{1/2}$.  For Lipschitz $h$
with $\omega_r\asymp k_r^{-1/s}$ on $\Xcal\subseteq\R^s$, balancing the
two gives the formal optimum $k_r\asymp n^{s/(s+2)}$ and rate
$O\{m\,n^{-1/(s+2)}\}$ up to logarithmic factors; for $s=1$,
$k_r\asymp n^{1/3}$ and rate $O(m\,n^{-1/3})$.  This is a heuristic
calibration, not a proved uniform theorem: the constants above are not
controlled uniformly in $r$ here, and the required uniform bounds are
listed as open in the Appendix.  In the simulation of
\S\,A11, the fitted exponent is $0.26$, below the $1/3$ heuristic,
while the total-error proxy decays as $n^{-0.336}$; the objective is
shallow near its discrete minimizer and the Monte Carlo standard deviation
of the induced conditional-likelihood estimator varies little over the
reported partition sizes, so we read these fits as descriptive.
\end{remark}

\begin{corollary}[Lifted diagonal regime]\label{cor:joint}
Fix the block size $m$ and lift the chart predictive law to $\Xcal^m$ by
\[
\widetilde\mu_{n,m}^{(r)}
=\mu_{n,m}^{(r)}\kappa_r^{\otimes m},
\]
where independently applying $\kappa_{r,j}=Q(\,\cdot\mid A_j^{(r)})$ fills
in each observed cell label.  For every fixed $r$, the window bound of
Theorem~\ref{thm:window} gives
$\TV\{\mu_{n,m}^{(r)},(p_r^*)^{\otimes m}\}\to0$ as $n\to\infty$, for any
tolerance with $\eta_n\to0$ and $n\eta_n\to\infty$, without lattice
regularity.  Under the projection approximation above, there is a
deterministic sequence $r_n\to\infty$ such that
\[
\TV\!\left\{\widetilde\mu_{n,m}^{(r_n)},(P^*)^{\otimes m}\right\}\to0 .
\]
Thus the finite-chart collapse recovers the ambient predictive law on the
original sample space, rather than only its cell probabilities.  This is a
projected-law conclusion; an experiment-level diagonal would require joint
control of the fixed-chart threshold and refinement error, which is not
available here.
\end{corollary}

\begin{proof}
For each chart, let $M_r$ be the least integer such that
$\TV(\mu_{n,m}^{(r)},(p_r^*)^{\otimes m})\le1/r$ for every $n\ge M_r$;
such an integer exists because the fixed-$r$ bound in
Theorem~\ref{thm:window} tends to zero.  To ensure that all preceding
thresholds have also been crossed, recursively set
$N_r=\max\{r,N_{r-1},M_r\}$, with $N_0=0$.  For
$n\ge N_1$ put $r_n=\max\{r:N_r\le n\}$, and set $r_n=1$ otherwise.
The maximum is finite because $N_r\ge r$, and $r_n\to\infty$.  Then
\begin{align*}
\TV\{\widetilde\mu_{n,m}^{(r_n)},(P^*)^{\otimes m}\}
&\le
\TV\{\mu_{n,m}^{(r_n)}\kappa_{r_n}^{\otimes m},
      (p_{r_n}^*)^{\otimes m}\kappa_{r_n}^{\otimes m}\}\\
&\quad+\TV\{(\widetilde P_{r_n}^*)^{\otimes m},(P^*)^{\otimes m}\}\\
&\le \tfrac1{r_n}+O(m\omega_{r_n})\to0,
\end{align*}
where the second inequality uses contraction of total variation under a
Markov kernel and Lemma~A3.
\end{proof}

\section{From Chart Curvature to Parameter Information}\label{sec:connections}

Whereas the diagonal argument concerns approximation across charts, the
statistical meaning of $H_r^*$ is a separate and local matter.  The
following score identity identifies what the reduced Hessian measures
before the geometry is pulled back to a parameter space.

\begin{remark}[Chart score experiment]\label{cor:lan}
If $J\sim p_r^*$, the categorical tangent score
$S_r(J)=V_r^\top\diag(1/p_r^*)e_J$ has mean zero and covariance
$\mathrm{Cov}\{S_r(J)\}=H_r^*$.  This exact finite identity explains why
the same quadratic form appears in Theorem~\ref{thm:gdf}; it does not
prove local asymptotic normality for the conditioned block law, which
remains open for growing blocks, including the finite-population
correction when $m/n$ does not vanish.
\end{remark}

In this limited sense the construction meets the local asymptotic theory
of \citet{lecam1972lan}: the reduced Hessian $H_r^*$ is the information
matrix for the tangent-coordinate experiment
$p_r^*+V_r v/\sqrt{n}$ and governs movement along the fixed constraint
surface, but is not in general the Fisher information for a parameter.

Parameter-space information enters only after this chart geometry is
pulled back through a parametric family.  In the setting of
\S\,\ref{sec:parametric}, write
$J_{\theta,r}=\partial p_{\theta,r}^*/\partial\theta^\top$ for the
Jacobian of the projection family; the parameter information under
correct specification is
\begin{equation}\label{eq:param-info}
I_{\theta,r}=J_{\theta,r}^\top\,\diag(1/p_{\theta,r}^*)\,J_{\theta,r},
\end{equation}
and $J_{\theta_0,r}$ must have full column rank for local identification.
A parameter-concentration rate is governed by
$\lambda_{\min}(I_{\theta_0,r})$, not by
$\lambda_{\min}(H_{\theta_0,r}^*)$.

\section{Parametric Comparisons}\label{sec:parametric}

\subsection{Parameter-Indexed Constraint Families}

The present paper departs from fixed-event theory at the parametric layer.
For a fixed $\theta$, moment conditioning again
produces an information projection; parametric inference begins only when
these projections are compared across $\theta$.  When the constraints
depend on $\theta\in\Theta\subset\R^p$, as in
$\Ebb_P\{h(X,\theta)\}=\alpha(\theta)$, the fixed constraint set
$\Ccal(\alpha)$ is replaced by a parameter-indexed family
$\{\Ccal_\theta(\alpha)\}$, and the projection and curvature become
functions of $\theta$, namely $P_\theta^*$, $p_{\theta,r}^*$ and
$H_{\theta,r}^*$, all computable on the finite chart.  The problem thus
acquires two layers: for
each $\theta$ one solves a constrained projection, then compares the
family across $\theta$ to infer the parameter.  The general form
$\Ebb_P\{g(X,\theta)\}=0$ is accommodated by $h_\theta(x)=g(x,\theta)$,
$\alpha(\theta)\equiv 0$ at the projection layer; additive separability
plays no role in the constraint geometry.  The fixed-matrix quadratic
expansion of \S\,\ref{sec:gmm} is narrower and treats the moving-moment
case separately.

\subsection{Induced Conditional-Likelihood Geometry}\label{sec:betel}

Let $\{p_{\theta,r}^*:\theta\in\Theta\}$ denote the projected family with
constraint target $\alpha(\theta)$, and let
$\mu_{n,m}^{(\theta,r)}$ denote the predictive law~\eqref{eq:predictive}
at that target.

\begin{definition}[Induced conditional likelihood]\label{def:icl}
The \emph{induced conditional likelihood} on the partition is the product
scoring rule associated with the limiting projected family,
\begin{equation}\label{eq:icl-def}
L_r^{\text{\textsc{icl}}}(\theta;\,x_{1:n})
=\prod_{i=1}^n p_{\theta,r}^*(T_r x_i),
\qquad T_r:\Xcal\to\{1,\ldots,k_r\}.
\end{equation}
It is a generalized predictive criterion induced by the limiting family of
Theorem~\ref{thm:collapse} or Theorem~\ref{thm:window}, as applicable; it
is not the likelihood of the conditioning event.
\end{definition}

\begin{proposition}[Predictive-risk calibration]\label{prop:plugin}
Fix a moment-adequate partition $\Pi_r$ and a block size $m$, and suppose
the constant of Theorem~\ref{thm:collapse}, or of Theorem~\ref{thm:window}
for a generic window, is bounded uniformly over a compact $K\subset\Theta$.
Then for every bounded
$\phi:\{1,\ldots,k_r\}^m\to\R$, uniformly over $\theta\in K$,
\begin{equation}\label{eq:plugin}
\bigl|\Ebb_{\mu_{n,m}^{(\theta,r)}}\phi
-\Ebb_{(p_{\theta,r}^*)^{\otimes m}}\phi\bigr|
\le 2\|\phi\|_\infty\,
\|\mu_{n,m}^{(\theta,r)}-(p_{\theta,r}^*)^{\otimes m}\|_{\TV},
\end{equation}
with constant one for $\phi\in[0,1]$.  The right-hand side is controlled by
Theorem~\ref{thm:collapse} in the exact-lattice case and by
Theorem~\ref{thm:window} for a generic window.
\end{proposition}

Thus the product prediction of $L_r^{\text{\textsc{icl}}}$ matches the
conditioned exchangeable-block law up to the collapse error.  In the
ambient setting the analogous product criterion is
$\prod_{i=1}^n(dP_\theta^*/dQ)(x_i)$.  The \emph{induced
conditional-likelihood estimator} is the maximizer over the feasible set,
\begin{equation}\label{eq:icl-est}
\widehat\theta_{\text{\textsc{icl}}}
=\argmax_{\theta:\,\Ccal_{\theta,r}(\alpha)\neq\varnothing}
\ \sum_{j=1}^{k_r} n_j\log p_{\theta,r,j}^*,
\qquad n_j=\#\{i:T_r X_i=j\},
\end{equation}
where $p_{\theta,r}^*$ is the Gibbs projection of the fixed baseline $q_r$
onto $\{p:A_r p=\alpha(\theta)\}$ of \S\,\ref{sec:geometry}; the estimator
uses the data only through the cell counts $(n_1,\ldots,n_{k_r})$.  Since
those counts are multinomial, \eqref{eq:icl-est} is the exact likelihood of
the cell counts under the induced family $\{p_{\theta,r}^*\}$, a smooth
positive submodel of the simplex.  For fixed $r$, a correctly specified
chart, and the identification and consistency conditions of
Proposition~A1, $\widehat\theta_{\text{\textsc{icl}}}$ is asymptotically normal
with variance $I_{\theta_0,r}^{-1}$ and is efficient \emph{within the chart
experiment} in the sense of the convolution theorem \citep{hajek1970}.
Efficiency in the ambient
experiment is not claimed: the discretization cost is the Le~Cam deficiency
of \S\,\ref{sec:lecam}, and the centring is at $\theta_0$ only when
$\sqrt n\,\omega_r\to0$, which the calibration of
Remark~\ref{rem:partition-sizing} does not satisfy.  Details are in the
Appendix, \S\,A6.2.

\begin{definition}[Induced generalized posterior]\label{thm:betel}
Given a prior $\pi$ on $\Theta$, the induced generalized posterior on the
partition is the exact update of $\pi$ by the induced conditional
likelihood of Definition~\ref{def:icl},
\begin{equation}\label{eq:betel}
\pi_r^{\text{\textsc{icl}}}(d\theta\mid X_{1:n})
=\frac{\pi(d\theta)\,L_r^{\text{\textsc{icl}}}(\theta;\,X_{1:n})\,
\ind\{\Ccal_{\theta,r}(\alpha)\neq\varnothing\}}
{\int\pi(d\vartheta)\,L_r^{\text{\textsc{icl}}}(\vartheta;\,X_{1:n})\,
\ind\{\Ccal_{\vartheta,r}(\alpha)\neq\varnothing\}} .
\end{equation}
The indicator restricts to $\theta$ for which the constrained projection
exists.  This is a generalized-Bayes construction, distinct in general
from the exact Bayesian posterior of the moment-conditioning event.
\end{definition}

Both $\widehat\theta_{\text{\textsc{icl}}}$ and its generalized posterior arise
directly from the generalized predictive criterion associated with the
collapse limit (Theorem~\ref{thm:collapse}).  On the finite-partition model the Gibbs
tilt $p_{\theta,r}^*$ shares the exponential-tilt form of the exponentially
tilted empirical likelihood of \citet{schennach2005}; the two coincide
exactly only in the separate empirical-support construction ($Q=\whatP$,
cells the data atoms), in which the $n$ empirical atoms are tilted with
baseline weights $1/n$.  For continuous data the connection is a geometric
analogy mediated by the partition, not an identity.  The feasibility
indicator in~\eqref{eq:betel} is the convex-hull condition of
\citet{schennach2005} and \citet{chibsimoni2022}.  Treating the tilted
empirical likelihood as a pseudo-likelihood,
\citet{chibsimoni2018,chibsimoni2022} establish the Bernstein--von~Mises
limit under much more general conditions.

Under standard regularity conditions for generalized-Bayes inference,
including prior mass near $\theta_0$, local risk separation, a uniform law
of large numbers, differentiability and stochastic equicontinuity, and
uniform lower bounds on the projected cell probabilities, the induced
generalized posterior of Definition~\ref{thm:betel} concentrates at
$\theta_0$ with local parameter information the pullback $I_{\theta_0,r}$
of~\eqref{eq:param-info}.  The relevant chart-to-parameter implication is
that $I_{\theta,r}$, not $\lambda_{\min}(H_r^*)$, governs parameter
concentration; contraction itself is standard, and the corresponding
asymptotics for the tilted empirical likelihood are given by
\citet{chibsimoni2018}.

\subsection{Generalized Method of Moments}\label{sec:gmm}

The classical generalized method of moments minimizes the quadratic
criterion
$\bar{g}_n(\theta)^\top W_n\bar{g}_n(\theta)$ where
$\bar{g}_n(\theta)=n^{-1}\sum_{i=1}^n h(X_i,\theta)-\alpha(\theta)$.
A quadratic criterion of the same form is induced by the local
projection geometry, as follows.  On the chart, write
$\bar{g}_{n,r}(\theta)=A_r\,\whatP^{(r)}-\alpha(\theta)$ for
the cell-level sample moments, where $A_r$ is the fixed cell-moment
matrix, and let
$\Sigma_{\theta,r}=\diag(p_{\theta,r}^*)
-p_{\theta,r}^*p_{\theta,r}^{*\top}$
denote the multinomial covariance at the projection, the inverse of
the curvature operator $\diag(1/p_{\theta,r}^*)$ on the simplex chart.

\begin{proposition}[Projection-induced moment criterion]\label{thm:gmm}
Fix a cell-moment matrix $A_r$ and a twice continuously differentiable target
$\alpha(\theta)$.  Suppose the cell labels are independent with law
$p_{\theta_0,r}^*$, $\theta_0$ is locally identified, and the
cell-moment covariance $\Omega_{\theta,r}=A_r\Sigma_{\theta,r}A_r^\top$ is
nonsingular near $\theta_0$.  Suppose also that, with probability tending
to one, $A_r\whatP^{(r)}$ lies in the relative interior of the finite mean
domain, so its dual solution exists.  Then uniformly over
$\|\theta-\theta_0\|_2\le Cn^{-1/2}$ and up to an additive term not
depending on $\theta$,
\begin{equation}\label{eq:gmm}
-\log L_r^{\text{\textsc{icl}}}(\theta)
=\frac{n}{2}\,\bar{g}_{n,r}(\theta)^\top
W^{\text{\textsc{icl}}}_{\theta,r}\,\bar{g}_{n,r}(\theta)
+O_{\Pbb}(n^{-1/2}),
\qquad
W^{\text{\textsc{icl}}}_{\theta,r}
=\{A_r\,\Sigma_{\theta,r}\,A_r^\top\}^{-1},
\end{equation}
so the induced weight is the inverse cell-moment covariance
$\Omega_{\theta,r}=A_r\Sigma_{\theta,r}A_r^\top$ under
the projected law.  Differentiating $A_r p_{\theta,r}^*=\alpha(\theta)$
gives the parameter information
$I_{\theta,r}=\dot\alpha(\theta)^\top\Omega_{\theta,r}^{-1}
\dot\alpha(\theta)$, so the local quadratic expansion of the induced
conditional likelihood carries the same optimal covariance weighting as a
generalized method-of-moments criterion for the projected cell moments.
\end{proposition}

The dual-expansion proof is in the Appendix, \S\,A6.3.

For a general parameter-dependent moment map $g(x,\theta)$ the matrix
$A_{\theta,r}$ itself varies with $\theta$ and the tangent space moves;
differentiating $A_{\theta,r}p_{\theta,r}^*=\alpha(\theta)$ then produces
an additional $\dot A_{\theta,r}p_{\theta,r}^*$ term.  Treating this case
would require uniform-in-$\theta$ control of the varying cell-moment
matrices and their derivatives, beyond the fixed-matrix expansion here.

On an empirical-support chart, with cells the observed support points and
$Q=\whatP$, $\Omega_{\theta,r}=A_r\Sigma_{\theta,r}A_r^\top$ becomes the
tilted empirical covariance of the moment function; under covariance
regularity, differentiability, local identification and stochastic
equicontinuity, \eqref{eq:gmm} then matches the efficient generalized method-of-moments criterion.  This
empirical-support construction is a separate, data-dependent
pseudo-likelihood chart (Remark~\ref{rem:baseline}), outside the
fixed-$Q$ theory.

\section{Parametric Misspecification}\label{sec:misspec}

When the data-generating law $P_0$ does not lie in the constraint
family, the target is the chart risk minimizer
$\theta_{0,r}=\argmin_{\theta\in\Theta}R_r(\theta)$ with
$R_r(\theta)=-\Ebb_{P_0}\log p_{\theta,r}^*(T_r X)$, a finite-chart
analogue of the Kullback--Leibler projection of \citet{Berk66} and
\citet{White1982}.  Two limits are involved.  For each fixed $\theta$, the
conditional block law collapses around $p_{\theta,r}^*$ by
Theorem~\ref{thm:collapse} or~\ref{thm:window}; separately, inference under
$P_0$ selects $\theta_{0,r}$.  Uniform collapse on a neighbourhood of
$\theta_{0,r}$, as in Proposition~\ref{prop:plugin}, therefore calibrates
the plug-in predictive law $(p_{\widehat\theta,r}^*)^{\otimes m}$ around
$(p_{\theta_{0,r},r}^*)^{\otimes m}$.  In the Gaussian illustration of
\S\,\ref{sec:simulation}, $\theta_0=0$ and $p_{\theta_0,r}^*=q_r$,
while $p_{M,r}^*$, the projection of the Model-M cell law $q_{M,r}$,
serves only to display the change of measure
(Fig.~\ref{fig:geometry-projection}).

Under a standard M-estimation and generalized-Bayes framework the induced
generalized posterior concentrates at $\theta_{0,r}$, its spread governed
by the parameter curvature $K_r=\nabla_\theta^2 R_r(\theta_{0,r})$ (assumed
positive definite) rather than by $H_{\theta_0,r}^*$, with the frequentist
estimator asymptotically normal with sandwich covariance
$n^{-1}K_r^{-1}\Sigma_r K_r^{-1}$,
$\Sigma_r=\mathrm{Var}_{P_0}\{\nabla_\theta\log p_{\theta_{0,r},r}^*(T_r X)\}$.
We do not develop this framework; the essential point is the distinction
between the parameter curvature $K_r$ and the tangent curvature $H_r^*$.

\section{Computational Illustration}\label{sec:simulation}

\subsection{Direct Check of Conditional Collapse}\label{sec:sim-direct}

\begin{figure}[H]
  \centering
  \includegraphics[width=0.52\textwidth]{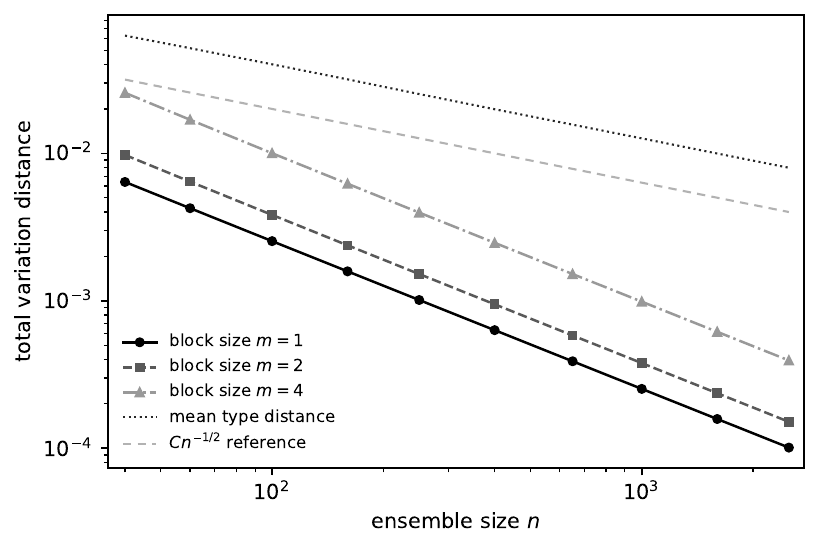}
  \caption{Exact conditional block collapse on a three-cell chart.
    Solid and broken curves give
    $\TV\{\mu_{n,m},(p^*)^{\otimes m}\}$ for fixed
    $m\in\{1,2,4\}$; the dotted curve is the conditional mean
    $\Ebb\{\TV(\widehat p,p^*)\}$, and the grey broken line is an
    $n^{-1/2}$ reference.  The type spread has the Gaussian scale from
    Theorem~\ref{thm:gdf}, while every fixed-block distance vanishes as
    asserted by Theorem~\ref{thm:collapse}.}
  \label{fig:direct-collapse}
  \alttext{On log-log axes, the conditional mean distance of the empirical
  type follows a square-root reference, while exact block-law distances
  for block sizes one, two and four decrease more rapidly.}
\end{figure}

We first isolate the mechanism of Theorem~\ref{thm:collapse} from the
parametric layer.  Consider three cells with baseline
$q=(0{.}2,0{.}5,0{.}3)$, cell moments $h=(-1,0,1)$ and target
$\alpha=0$.  Exact feasibility is the integer constraint $N_1=N_3$, and
the information projection is
$p^*=(0{.}24745,0{.}50510,0{.}24745)$.  For each $n$ we enumerate every
feasible type $(t,n-2t,t)$, evaluate its conditional multinomial weight
in~\eqref{eq:typeweight}, and mix the corresponding ordered
hypergeometric laws in~\eqref{eq:mixture}.  Thus the distances in
Fig.~\ref{fig:direct-collapse} contain no Monte Carlo error.

The mean type distance follows the expected $n^{-1/2}$ scale.  The block
distances decrease faster in this example because mixing cancels the
first-order centred type fluctuation; Theorem~\ref{thm:collapse} supplies
an upper bound rather than an exact asymptotic rate.  Enumeration details
are in the Appendix, \S\,A9.

\subsection{Parametric Design}

The second illustration is designed to display the parametric geometry
rather than to rank the estimators.  We adapt the mean-and-variance
constraint problem
of \citet{schennach2007el}, Section~4, retaining its moment functions but
using a more pronounced misspecified regime.  For $X_1,\ldots,X_n$ from a
location--scale normal model, the moment functions
$g(x,\theta)=\{x-\theta,\,(x-\theta)^2-1\}$ give a constraint
$E\{g(X,\theta_0)\}=0$ that identifies the mean $\theta_0$ while pinning
the variance to one.  We consider two regimes: Model~C, with
$X_i\sim N(0,1)$, correctly specified; and Model~M, with $X_i\sim N(0,4)$,
in which the variance constraint is misspecified (a fourfold variance
ratio).  The baseline $Q$ is $N(0,1)$ and we use a quantile partition of
size $k_r=20$, on which $q_r$ is uniform.  Because the extreme cells are
unbounded, Assumption~\ref{ass:refine} does not apply to this fixed chart.
We therefore treat $k_r=20$ as a fixed-chart calculation and compute the
finite cell moments $h_r(j)$ in closed form from the truncated-normal
formulae.

\begin{samepage}
The four laws used below have distinct roles:
\begin{center}
\small
\begin{tabular}{@{}ll@{}}
$q_r$ & fixed baseline cell law; also $p_{\theta_0,r}^*$ in this design,\\
$q_{M,r}$ & cell law of the misspecified data-generating distribution,\\
$p_{M,r}^*$ & projection of $q_{M,r}$, used for change-of-measure geometry,\\
$p_{\theta,r}^*$ & parameter-indexed projection used for inference.
\end{tabular}
\end{center}
\end{samepage}

The comparison adapts the mean--variance design of \citet{schennach2007el}.
The two tilted procedures enforce both moments but differ in support: the
ETEL maximizer of the BETEL profile \citep{schennach2005}, labelled
$\hat\theta_{\text{\textsc{betel}}}$, tilts the realized sample, whereas
$\hat\theta_{\text{\textsc{icl}}}$ tilts the fixed partition cells.  The
separate benchmark $\hat\theta_{\text{\textsc{gmm}}}=\bar X$ uses only the
correctly specified location moment and leaves the variance restriction
unused.  This variance-agnostic reference is specific to the present paper,
not to \citet{schennach2007el} (Appendix, \S\,A9).

\subsection{Geometric Structure}\label{sec:sim-geometry}

\begin{figure}[htbp]
  \centering
  \includegraphics[width=0.55\textwidth]{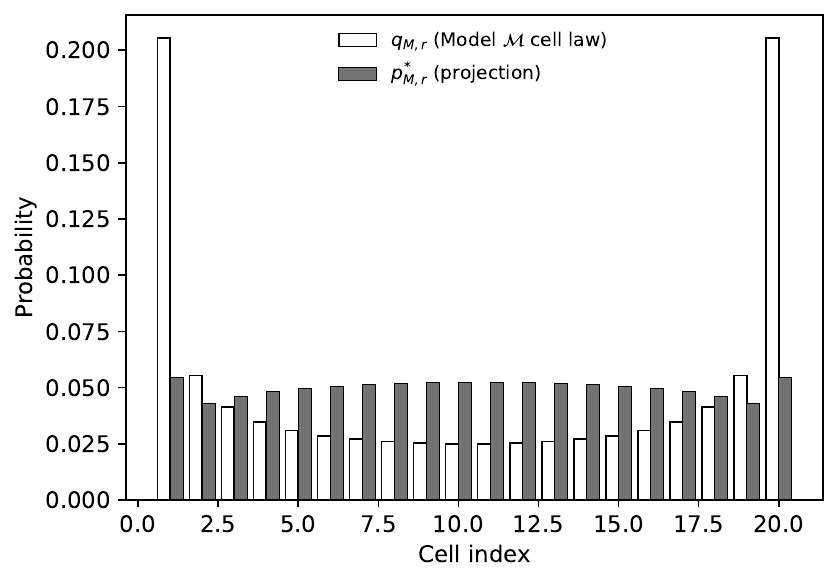}
  \caption{Kullback--Leibler projection under Model~M ($k_r=20$,
    $\sigma=2$): cell law
    $q_{M,r}$ of the $N(0,4)$ measure and its information projection
    $p_{M,r}^*$, which redistributes the extreme-cell mass towards the
    centre to satisfy the variance constraint ($\KL\approx 0{.}269$).}
  \label{fig:geometry-projection}
  \alttext{Paired bars show the misspecified law placing more mass in extreme
  cells and its projection shifting mass towards the centre.}
\end{figure}

Figure~\ref{fig:geometry-projection} shows the Model-M cell law
$q_{M,r}$ and its information projection $p_{M,r}^*$ onto the constraint
set.  This change-of-measure projection is distinct from the parametric
pseudo-true law $p_{\theta_0,r}^*$ of \S\,\ref{sec:misspec}; here the
baseline cell law already satisfies the constraints, so
$p_{\theta_0,r}^*=p_r^*=q_r$.  The variance restriction moves mass from
the extreme cells towards the centre, with
$\KL(p_{M,r}^*\|q_{M,r})\approx0{.}269$.

\begin{figure}[htbp]
  \centering
  \begin{minipage}[t]{0.45\textwidth}
    \centering
    \includegraphics[width=\linewidth]{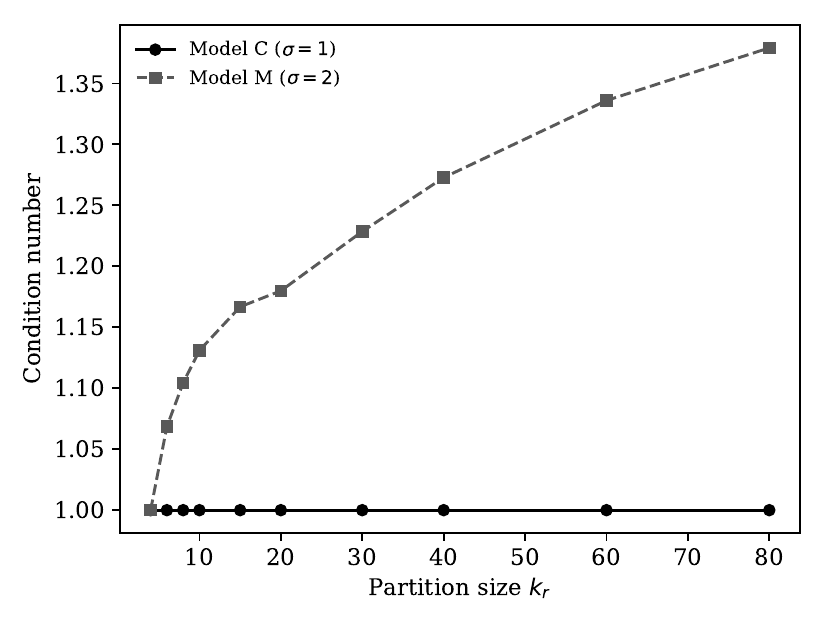}
  \end{minipage}\hfill
  \begin{minipage}[t]{0.45\textwidth}
    \centering
    \includegraphics[width=\linewidth]{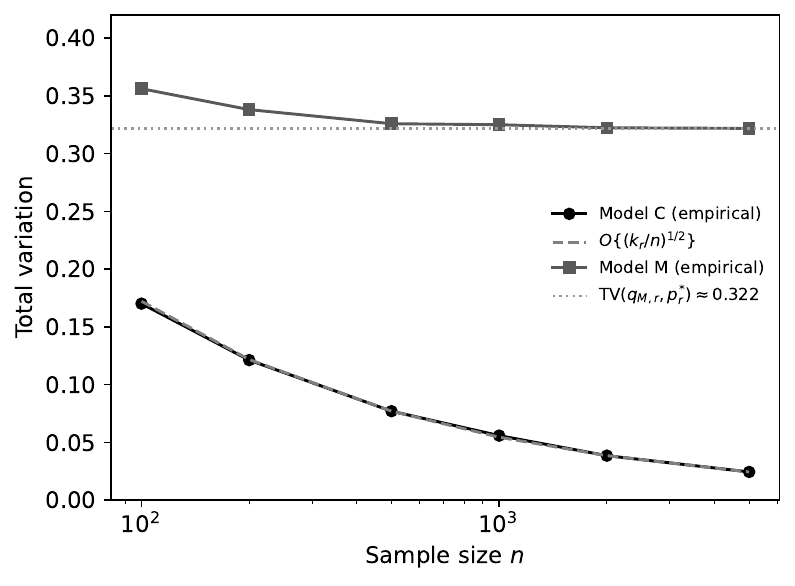}
  \end{minipage}
  \caption{Left: conditioning $\kappa(H_{M,r}^*)$ of the Model-M
    change-of-measure projection $q_{M,r}\to p_{M,r}^*$ ($H_{M,r}^*$ its
    reduced Hessian, distinct from the baseline $H_r^*$), isotropic under
    Model~C and growing to $\kappa\approx 1{.}379$ at $k_r=80$ under
    Model~M.  Right: empirical
    cell discrepancy $\TV(\widehat p_n,p_r^*)$ versus $n$ ($k_r=20$,
    $500$ replications); under Model~C it tracks $O\{(k_r/n)^{1/2}\}$,
    under Model~M it converges to the floor
    $\TV(q_{M,r},p_r^*)\approx 0{.}322$.}
  \label{fig:curv-collapse}
  \alttext{Left, the condition number stays flat under Model C and rises under
  Model M.  Right, total variation vanishes under Model C but levels near
  0.322 under Model M.}
\end{figure}

Figure~\ref{fig:curv-collapse} gives two complementary diagnostics.
The condition number of the change-of-measure curvature is one under
Model~C and reaches $1{.}379$ at $k_r=80$ under Model~M.  The mean
empirical discrepancy $\TV(\widehat p_n,p_r^*)$ decays at the rate
$O\{(k_r/n)^{1/2}\}$ under Model~C but approaches
$\TV(q_{M,r},p_r^*)\approx0{.}322$ under Model~M.  This diagnostic
concerns misspecification and is distinct from the conditional
block-predictive law of Theorem~\ref{thm:collapse}; likewise, $p_{M,r}^*$
is a separate object, with $\TV(q_{M,r},p_{M,r}^*)\approx0{.}326$.

\subsection{Comparison of Estimation Procedures}\label{sec:comparison-est}

Table~\ref{tab:sim} reports the Monte Carlo standard deviation and
bias of each estimator over $2000$ replications at sample sizes
$n\in\{1000,5000\}$.  Under Model~C all three procedures behave
similarly, as expected when the moment conditions are correctly
specified.  Under Model~M a dispersion difference emerges:
$\hat\theta_{\text{\textsc{betel}}}$ becomes more dispersed, while
$\hat\theta_{\text{\textsc{icl}}}$ remains close to $\hat\theta_{\text{\textsc{gmm}}}$
(Appendix, Fig.~S1).  Biases are
negligible in both regimes.

\begin{table}[htbp]
  \centering
  \caption{Monte Carlo standard deviation ($\times 10^2$) and bias
    ($\times 10^4$) for estimation of $\theta_0=0$.
    Model~C: $X_i\sim N(0,1)$; Model~M: $X_i\sim N(0,4)$.
    Each entry is based on $2000$ replications; Monte Carlo standard
    errors of the bias entries range from $3$ to $23$ in the same
    units.  The final column is the mean total-variation discrepancy
    $\TV(\widehat p_n,p_r^*)$ of \S\,\ref{sec:sim-geometry}.
    Computational details are
    in the Appendix, \S\,A9}\label{tab:sim}
  \begin{tabular}{llrrrrrrr}
    \toprule
    && \multicolumn{3}{c}{Standard deviation} &
    \multicolumn{3}{c}{Bias} & Mean \\
    \cmidrule(lr){3-5}\cmidrule(lr){6-8}
    Model & $n$ & \textsc{gmm} & \textsc{betel} & \textsc{icl} &
    \textsc{gmm} & \textsc{betel} & \textsc{icl} & TV \\
    \midrule
    C & 1000 & 3.157 & 3.170 & 3.122 & $-5.7$ & $-6.1$ & $-5.8$ & 0.055 \\
    C & 5000 & 1.406 & 1.407 & 1.402 & $2.0$ & $2.0$ & $2.4$ & 0.025 \\
    M & 1000 & 6.409 & 10.329 & 6.231 & $-7.4$ & $-19.8$ & $-6.3$ & 0.324 \\
    M & 5000 & 2.792 & 4.539 & 2.711 & $3.5$ & $2.6$ & $5.6$ & 0.322 \\
    \bottomrule
  \end{tabular}
\end{table}

Under Model~M, the last column of Table~\ref{tab:sim} approaches the floor
$0{.}322$.  This comparison illustrates the distinction between
specification and misspecification; the relative performance of the
estimators remains specific to this design.

\FloatBarrier

\section{Discussion}\label{sec:discussion}

The central gain from conditioning is a probabilistic interpretation of the
information projection.  The exact type mixture separates two mechanisms:
the conditional empirical type concentrates near $p_r^*$, and sampling a
fixed block without replacement becomes sampling from its cell
probabilities.  Lattice regularity sharpens the first mechanism to a
Gaussian statement; the shrinking-window argument retains the predictive
conclusion when that arithmetic structure is absent.

The reference law remains indispensable.  Moment restrictions alone do not
identify a predictive distribution, and changing $Q$ changes both the
projection and its curvature.  Refinement controls the chart approximation
for a fixed reference construction; it does not remove this modelling
choice.  Parameterization is a second operation: differentiating the
projected family produces parameter information and the local
generalized-method-of-moments criterion.  This explains why tangent
curvature, parameter information and misspecification risk curvature play
different roles.

The main limitation is that the concentration results are proved on a fixed
chart, with constants that need not remain controlled under refinement.  A
uniform theory in which chart dimension and block size grow together would
require stable normalization of the curvature spectrum and joint control of
refinement, localization and finite-population error.  Such control would
replace the present diagonal existence argument by an explicit joint
regime.  Partial exchangeability and sequential prediction provide natural
further settings for the same probabilistic construction.

\section*{Acknowledgement}
This work originated from a 2012 seminar at the University of Texas at
Austin.

\clearpage

\setcounter{section}{0}
\setcounter{theorem}{0}
\setcounter{lemma}{0}
\setcounter{corollary}{0}
\setcounter{proposition}{0}
\setcounter{definition}{0}
\setcounter{assumption}{0}
\setcounter{remark}{0}
\setcounter{equation}{0}
\setcounter{table}{0}
\setcounter{figure}{0}
\renewcommand{\thesection}{A\arabic{section}}
\renewcommand{\thetheorem}{A\arabic{theorem}}
\renewcommand{\thelemma}{A\arabic{lemma}}
\renewcommand{\thecorollary}{A\arabic{corollary}}
\renewcommand{\theproposition}{A\arabic{proposition}}
\renewcommand{\thedefinition}{A\arabic{definition}}
\renewcommand{\theassumption}{A\arabic{assumption}}
\renewcommand{\theremark}{A\arabic{remark}}
\renewcommand{\theequation}{A\arabic{equation}}
\renewcommand{\thetable}{A\arabic{table}}
\renewcommand{\thefigure}{A\arabic{figure}}

\begin{center}
{\large\bfseries Appendix}
\end{center}
\medskip

This supplement gives the proofs, technical lemmas and simulation details
supporting the main paper.  It also develops the Le~Cam comparison in full.
Supplementary results carry the prefix~S; unprefixed theorem, proposition,
definition and equation numbers refer to the main paper.

\section*{Notation}\label{sec:notation}

Symbols are grouped by the layer of the construction in which they first
appear.  References of the form \S\,2.1, Definition~1, Theorem~2 or~(13)
are to the \emph{main paper}; references of the form \S\,A4 are to this appendix.  A subscript $r$ always indexes the partition (chart), $n$ the
sample size and $m$ the block size.

\begingroup
\footnotesize
\setlength{\LTpre}{4pt}\setlength{\LTpost}{4pt}
\begin{longtable}{@{}l p{8.0cm} l@{}}
\toprule
Symbol & Meaning & Defined at \\
\midrule
\endfirsthead
\toprule
Symbol & Meaning & Defined at \\
\midrule
\endhead

\multicolumn{3}{@{}l}{\itshape Ambient objects}\\[2pt]
$(\Xcal,\Acal)$ & sample space and its $\sigma$-field & \S\,2.1 \\
$Q$ & baseline probability measure on $\Xcal$ & \S\,2.1 \\
$h:\Xcal\to\R^d$ & moment map; $d$ the number of moments & \S\,2.1 \\
$\alpha\in\R^d$ & moment target & \S\,2.1 \\
$\Ccal(\alpha)$ & admissible class $\{P\ll Q:\int h\,dP=\alpha\}$ & (1) \\
$P^*$ & information projection of $Q$ onto $\Ccal(\alpha)$ & Thm.~1 \\
$\lambda^*,\ \psi$ & ambient tilt and log-partition function & Thm.~1, (2) \\
$\KL(\cdot\|\cdot)$ & Kullback--Leibler divergence & \S\,2.1 \\
\addlinespace

\multicolumn{3}{@{}l}{\itshape Partition chart}\\[2pt]
$\Pi_r$ & finite measurable partition (chart) of $\Xcal$ & \S\,2.3 \\
$A_j^{(r)},\ k_r$ & cells of $\Pi_r$ and their number & \S\,2.3 \\
$T_r$ & partition map $\Xcal\to\{1,\ldots,k_r\}$ & \S\,3.1 \\
$q_r,\ q_{r,j}$ & baseline cell law, $q_{r,j}=Q(A_j^{(r)})$ & Def.~1 \\
$h_r(j)$ & cell moment $\Ebb_Q\{h(X)\mid X\in A_j^{(r)}\}$ & Def.~1 \\
$A_r$ & cell-moment matrix $[h_r(1)\cdots h_r(k_r)]$, $d\times k_r$ & Def.~1 \\
$\widetilde A_r$ & augmented constraint matrix $(\ind\;A_r^\top)^\top$ & Def.~1 \\
$\Ccal_r(\alpha)$ & chart constraint set $\{p:A_r p=\alpha\}$ & \S\,2.3 \\
$p_r^*,\ p_{\min,r}^*$ & projected cell law and its smallest cell probability & \S\,2.3--2.4 \\
$\lambda_r^*,\ \psi_r$ & chart tilt and its log-partition function & \S\,2.3, \S\,A1 \\
\addlinespace

\multicolumn{3}{@{}l}{\itshape Constraint geometry and curvature}\\[2pt]
$T_r^*$ & tangent space to the constraint surface at $p_r^*$ & \S\,2.4 \\
$V_r$ & $k_r\times s_r$ orthonormal basis of $T_r^*$ & \S\,2.4 \\
$s_r=k_r-1-d$ & dimension of $T_r^*$ & \S\,2.4, \S\,3.1 \\
$H_r^*$ & reduced Hessian $V_r^\top\diag(1/p_r^*)V_r$ & \S\,2.4 \\
$\lambda_{\min}(H_r^*),\ \kappa(H_r^*)$ & minimal curvature; condition number & \S\,2.4--2.5 \\
$\Sigma_{\theta,r}$ & multinomial covariance $\diag(p^*_{\theta,r})-p^*_{\theta,r}p_{\theta,r}^{*\top}$ & \S\,6.3 \\
$\Omega_{\theta,r}$ & $d\times d$ cell-moment covariance $A_r\Sigma_{\theta,r}A_r^\top$ & Prop.~3 \\
$J_{\theta,r},\ I_{\theta,r}$ & Jacobian of the projected family; parameter information & (18) \\
\addlinespace

\multicolumn{3}{@{}l}{\itshape Predictive representation}\\[2pt]
$Z_i=T_r(X_i)$ & cell label of observation $i$ & \S\,3.1 \\
$\whatP^{(r)}$ & empirical type $n^{-1}\sum_i e_{Z_i}$ & \S\,3.1 \\
$\Pcal_n$ & lattice of realizable $n$-types & \S\,3.1 \\
$E_{n,r}$ & feasible types $\{p\in\Pcal_n:\|A_r p-\alpha\|\le\eta_n\}$ & (6) \\
$\eta_n$ & feasibility tolerance ($\eta_n=0$: exact feasibility) & \S\,3.1 \\
$\mathcal{N}_r$ & cofinal set of exactly feasible $n$ (lattice regularity) & \S\,3.1 \\
$m$ & block size, $1\le m\le n$ & \S\,3.1 \\
$\mu_{n,m}^{(r)}$ & law of a selected exchangeable block & (7) \\
$w_{n,r}(p)$ & constrained multinomial type weight & (8) \\
$\mathrm{HG}_{n,p,m}$ & within-type hypergeometric law & (9) \\
$\bar p_{n,r}$ & \emph{discrete} minimizer over exactly feasible $n$-types & \S\,3.1, \S\,A4 \\
$p_r^{\eta_n}$ & \emph{continuous} minimizer over the thickened slab & Rem.~4, \S\,A4 \\
\addlinespace

\multicolumn{3}{@{}l}{\itshape Localization}\\[2pt]
$\mathcal{G}_{n,r},\ \mathcal{B}_{n,r}$ & good and bad types & (12) \\
$\rho_n$ & localization radius $\{4\log n/\lambda_{\min}(H_r^*)\}^{1/2}$ & \S\,3.2 \\
$\Lambda_r$ & tangent lattice of feasible types, mesh $n^{-1/2}$ & \S\,A3, Lem.~\ref{lem:riemann} \\
$R_{n,m,r}$ & remainder in the Gaussian localization & Thm.~2, (13) \\
$C_r$ & generic constant depending only on the chart $r$ & Thms.~2--4 \\
$\omega_r$ & cell oscillation $\sup_j\sup_{x,x'\in A_j^{(r)}}\|h(x)-h(x')\|$ & Ass.~2 \\
\addlinespace

\multicolumn{3}{@{}l}{\itshape Le~Cam comparison and refinement}\\[2pt]
$\Ecal_n,\ \Ecal_{n,r}$ & ambient and discretized experiments & \S\,4 \\
$\delta,\ \Delta$ & deficiency and Le~Cam distance & \S\,4, (17) \\
$\kappa_{r,j}$ & reverse kernel $Q(\cdot\mid A_j^{(r)})$ & \S\,4 \\
$\widetilde P_{\theta,r},\ \widetilde P_r^*$ & lifted law and lifted projection & \S\,4, \S\,A7 \\
$\widetilde\mu_{n,m}^{(r)}$ & chart predictive law lifted to $\Xcal^m$ & Cor.~1 \\
\addlinespace

\multicolumn{3}{@{}l}{\itshape Parametric layer}\\[2pt]
$\theta\in\Theta,\ \alpha(\theta)$ & parameter and parameter-indexed target & \S\,6.1 \\
$p_{\theta,r}^*,\ \lambda_{\theta,r}$ & projected family on the chart and its tilt & \S\,6.1, \S\,A6.2--A6.3 \\
$L_r^{\text{\textsc{icl}}}$ & induced conditional likelihood & Def.~2, (19) \\
$\widehat\theta_{\text{\textsc{icl}}}$ & induced conditional-likelihood estimator & (21) \\
$\pi_r^{\text{\textsc{icl}}}$ & induced generalized posterior & Def.~3, (22) \\
$n_j$ & cell counts $\#\{i:T_rX_i=j\}$ & \S\,6.2 \\
$\bar g_{n,r}(\theta)$ & chart sample moments $A_r\whatP^{(r)}-\alpha(\theta)$ & \S\,6.3 \\
$W^{\text{\textsc{icl}}}_{\theta,r}$ & induced quadratic-criterion weight $\Omega_{\theta,r}^{-1}$ & (23) \\
\addlinespace

\multicolumn{3}{@{}l}{\itshape Misspecification}\\[2pt]
$P_0$ & data-generating law & \S\,7 \\
$R_r(\theta),\ \theta_{0,r}$ & chart risk and its minimizer & \S\,7 \\
$K_r$ & parameter curvature $\nabla_\theta^2R_r(\theta_{0,r})$ & \S\,7 \\
$\Sigma_r$ & score variance in the sandwich covariance & \S\,7 \\
\addlinespace

\multicolumn{3}{@{}l}{\itshape Computational illustration}\\[2pt]
$g(x,\theta)$ & moment function $\{x-\theta,\,(x-\theta)^2-1\}$ & \S\,8.2 \\
$q_{M,r},\ p_{M,r}^*$ & Model-M cell law and its change-of-measure projection & \S\,8.2--8.3 \\
$H_{M,r}^*$ & reduced Hessian of the Model-M projection & \S\,8.3, Fig.~4 \\
\bottomrule
\end{longtable}
\addtocounter{table}{-1}
\endgroup
\renewcommand{\theHtable}{S\arabic{table}}

\section{Proof of Theorem~1 (Measure-Level Projection)}%
\label{sec:proof-projection}

The main paper states the projection result and its Pythagorean
interpretation; here we give the complete proof.

Under Assumption~1 the log-moment-generating function $\psi$ is
essentially smooth and steep, and
$\alpha\in\nabla\psi(\Theta_\psi)$.  The concave dual
$\lambda\mapsto\lambda^\top\alpha-\psi(\lambda)$ therefore attains its
maximum at the unique interior point $\lambda^*$ satisfying
$\nabla\psi(\lambda^*)=\alpha$; essential smoothness rules out a boundary
maximizer, while the negative-definite Hessian $-\nabla^2\psi(\lambda^*)$
gives uniqueness.  For unbounded $h$, we do not assume that
$\Ccal(\alpha)$ is weakly closed.  Existence and optimality instead follow
from duality and from the Pythagorean identity below.  We state the identity
first for feasible $P$ with finite $\KL(P\|Q)$; if
$\KL(P\|Q)=\infty$, the minimizing conclusion is immediate.  The tilt
$e^{\lambda^\top h(x)-\psi(\lambda)}Q(dx)$ is a probability measure for
$\lambda\in\Theta_\psi$, and the first-order (Gibbs) form
$\log(dP^*/dQ)=\lambda^{*\top}h-\psi(\lambda^*)$ follows.

To verify that $P^*\in\Ccal(\alpha)$, compute
$\int h\,dP^*=\int h(x)e^{\lambda^{*\top}h(x)-\psi(\lambda^*)}Q(dx)
=\nabla\psi(\lambda^*)=\alpha$
by definition of $\lambda^*$.  The value
$\KL(P^*\|Q)=\lambda^{*\top}\alpha-\psi(\lambda^*)$ follows by
substitution of the Gibbs density into the definition of the
Kullback--Leibler divergence.

We finally compare $P^*$ with an arbitrary feasible law.  For any
$P\in\Ccal(\alpha)$ with $P\neq P^*$ and $\KL(P\|Q)<\infty$,
\begin{align}
\KL(P\|Q)-\KL(P^*\|Q)
&=\int\log\frac{dP}{dQ}\,dP-\int\log\frac{dP^*}{dQ}\,dP^*\notag\\
&=\KL(P\|P^*)+\int\log\frac{dP^*}{dQ}(dP-dP^*)\notag\\
&=\KL(P\|P^*)+\lambda^{*\top}\left(\int h\,dP-\int h\,dP^*\right)
  \notag\\
&=\KL(P\|P^*)>0,\label{eq:strict}
\end{align}
where the penultimate equality uses $\int h\,dP=\int h\,dP^*=\alpha$.
This establishes global optimality.\qed

\section{Proof of Proposition~1 (Quadratic Expansion)}%
\label{sec:proof-quadratic}

The main paper records the first-order argument; here we give the full
third-order remainder calculation.

Write $\KL(p\|q_r)=\sum_j p_j\log(p_j/q_{r,j})$ and expand around
$p_r^*$ with $p=p_r^*+V_r u$.  The gradient at $p_r^*$ restricted to
$T_r^*$ vanishes by the KKT conditions, as shown in the main paper.
The Hessian of $\KL$ at $p$ is $\diag(1/p)$, which restricted to
$T_r^*$ and evaluated at $p_r^*$ gives $H_r^*$.

For the remainder, the third derivative of
$f(p_j)=p_j\log(p_j/q_{r,j})$ is $f'''(p_j)=-1/p_j^2$.  Under the
radius condition $\|u\|_2\le p_{\min,r}^*/2$ of the statement, every
coordinate along the segment $[p_r^*,p]$ satisfies
$\tilde{p}_j\ge p_{r,j}^*-\|V_ru\|_\infty\ge p_{\min,r}^*/2$, so by
the multivariate Taylor remainder theorem,
\[
|R_3(u)|\le\frac{1}{6}\sup_{\tilde{p}\in[p_r^*,p]}
\sum_{j=1}^{k_r}\frac{|(V_r u)_j|^3}{\tilde{p}_j^2}
\le\frac{4\,\|V_r u\|_3^3}{6(p_{\min,r}^*)^2}
\le\frac{2\,\|u\|_2^3}{3(p_{\min,r}^*)^2},
\]
where the last inequality uses $\|V_r u\|_3\le\|V_r u\|_2\le\|u\|_2$
since $V_r$ has orthonormal columns.\qed

\section{Proofs of Lemmas~1 and~2 (Localization)}%
\label{sec:proof-localization}

The geometry of the finite chart comes from the constraints in
\S\,2.4.  The localization window, used in the spirit of
\citet{Lanford1973}, serves a narrower purpose: it separates the region
where a quadratic approximation is accurate from the region where a
large-deviation bound suffices.

\subsection{Proof of Lemma~1 (Exponential Suppression)}

In the exact-feasibility case ($\eta_n=0$) the continuous constrained
minimizer over $\Ccal_r(\alpha)$ is $p_r^*$, whose Gibbs coordinates are
generally irrational, so $p_r^*$ need not be an $n$-type; under lattice
regularity (\S\,A4) the discrete minimizer over exactly feasible types
satisfies $\|\bar p_{n,r}-p_r^*\|_2=O_r(n^{-1})$.  We centre the analysis at
the continuous $p_r^*$; the $O_r(n^{-1})$ gap is negligible on the annulus
scale $\rho_n/\sqrt{n}$.  Split $\mathcal{B}_{n,r}$ into a local annulus
$\rho_n/\sqrt{n}<\|\widehat p-p_r^*\|_2\le\delta_r$ and a far region
$\|\widehat p-p_r^*\|_2>\delta_r$, where $\delta_r>0$ is fixed small enough
that Proposition~1 holds on the ball of radius $\delta_r$ and every cell
probability there exceeds $p_{\min,r}^*/2$.

By Proposition~1, decrease $\delta_r$ if necessary so that the cubic
remainder is at most one quarter of the quadratic term.  On the annulus,
\[
\KL(\widehat p\|q_r)-\KL(p_r^*\|q_r)\ge
\tfrac14\lambda_{\min}(H_r^*)\|\widehat p-p_r^*\|_2^2.
\]
With $\rho_n^2=4\log n/\lambda_{\min}(H_r^*)$ the excess is therefore at
least $\log n/n$ at the inner radius and grows outward.  On the far region, strict
convexity of $\KL(\cdot\|q_r)$ and uniqueness of the minimizer give a
fixed gap $\gamma_r(\delta_r)>0$.

The conditional bad mass is the ratio of summed multinomial weights over
$\mathcal{B}_{n,r}$ to those over $\mathcal{G}_{n,r}$.  We bound it by a
lattice-sum comparison, using Lemma~2 \emph{only} for the denominator.
Under lattice regularity the exactly feasible types form a translate of a
rank-$s_r$ lattice of mesh $O(n^{-1})$ in $p$, hence, in the scaled tangent
coordinate $u=\sqrt{n}\,V_r^\top(\widehat p-p_r^*)$, a lattice $\Lambda_r$
of mesh $O(n^{-1/2})$ and covolume $\asymp n^{-s_r/2}$.  Write
$w_r=(\prod_j p_{r,j}^*)^{-1/2}(2\pi n)^{-(k_r-1)/2}
\exp\{-n\KL(p_r^*\|q_r)\}$ for the common Stirling factor.

\medskip\noindent\textit{Denominator.}
For $\|u\|_2\le\rho_n$, Lemma~2 gives the good-type weight
$w_r\exp\{-\tfrac12 u^\top H_r^*u\}\{1+o(1)\}$, so by the lattice
Riemann-sum estimate of Lemma~A4 (\S\,A8),
\[
\sum_{u\in\Lambda_r,\ \|u\|_2\le\rho_n}
\exp\{-\tfrac12 u^\top H_r^* u\}
\ \asymp\ n^{s_r/2}\int_{\R^{s_r}}\exp\{-\tfrac12 u^\top H_r^* u\}\,du
\ \asymp\ n^{s_r/2}.
\]

\medskip\noindent\textit{Numerator.}
Lemma~2 is not used on the annulus.  There every cell probability exceeds
$p_{\min,r}^*/2$, so Stirling's formula gives the \emph{uniform} upper bound
$\Pbb(\whatP^{(r)}=\widehat p)\le C_r\,(\prod_j p_{r,j}^*)^{-1/2}
(2\pi n)^{-(k_r-1)/2}\exp\{-n\KL(\widehat p\|q_r)\}$, while the quadratic
lower bound converts
$\exp\{-n[\KL(\widehat p\|q_r)-\KL(p_r^*\|q_r)]\}$ into
$\exp\{-\tfrac14\lambda_{\min}(H_r^*)\|u\|_2^2\}$.  The annulus sum is
therefore at most
\[
C_r\,w_r\!\!\sum_{u\in\Lambda_r,\ \|u\|_2>\rho_n}\!\!
\exp\{-\tfrac14\lambda_{\min}(H_r^*)\|u\|_2^2\}
\ \le\ C_r\,w_r\,n^{s_r/2}\;\Pbb(\|Z\|_2>\rho_n),
\]
where $Z\sim\Nor\{0,2\lambda_{\min}(H_r^*)^{-1}I\}$,
again by Lemma~A4.  Both numerator and denominator contain
$w_r\,n^{s_r/2}$, so this factor drops out.  What remains is the Gaussian
tail
$\Pbb(\|Z\|_2>\rho_n)\le
C_r(1+\rho_n)^{s_r}e^{-\lambda_{\min}(H_r^*)\rho_n^2/4}
=O_r\{n^{-1}(\log n)^{s_r/2}\}$.  The far region
contributes at most $|\mathcal{B}_{n,r}|\,e^{-n\gamma_r(\delta_r)}$ relative
to the denominator, negligible for fixed $r$.  Hence the bad mass is
$O_r\{n^{-1}(\log n)^{s_r/2}\}$.\qed

\subsection{Proof of Lemma~2 (Gaussian Control)}

The multinomial probability of type $\widehat{p}$ is
\[
\Pbb(\whatP^{(r)}=\widehat{p})
=\binom{n}{n\widehat{p}_1,\ldots,n\widehat{p}_{k_r}}
\prod_{j=1}^{k_r}q_{r,j}^{n\widehat{p}_j}.
\]
Stirling's approximation for the multinomial coefficient gives
\[
\binom{n}{n\widehat{p}}\approx(2\pi n)^{-(k_r-1)/2}
\left(\prod_{j=1}^{k_r}\widehat{p}_j\right)^{-1/2}
\exp\{n\,H(\widehat{p})\},
\]
where $H(\widehat{p})=-\sum_j\widehat{p}_j\log\widehat{p}_j$ is the
Shannon entropy.  Combining with the baseline term gives
\[
\Pbb(\whatP^{(r)}=\widehat{p})
\propto\left(\prod_j\widehat{p}_j\right)^{-1/2}
\exp\{-n\,\KL(\widehat{p}\|q_r)\}.
\]

The Stirling prefactor $(\prod_j\widehat{p}_j)^{-1/2}$ at
$\widehat{p}=p_r^*+V_r u/\sqrt{n}$ satisfies, for fixed $r$,
\[
\log\frac{\prod_j\widehat p_j}{\prod_j p_{r,j}^*}
=\sum_j\log\Bigl(1+\frac{(V_r u)_j}{\sqrt n\,p_{r,j}^*}\Bigr)
=O_r\!\Bigl(\frac{\|u\|_2}{\sqrt n}\Bigr),
\]
the constant depending on $k_r$ and $p_{\min,r}^*$ through
$\sum_j|(V_r u)_j|/p_{r,j}^*\le\sqrt{k_r}\,\|u\|_2/p_{\min,r}^*$; hence,
for $\|u\|_2\le\rho_n$,
$(\prod_j\widehat p_j)^{-1/2}=(\prod_j p_{r,j}^*)^{-1/2}\{1+O_r(\rho_n/\sqrt{n})\}$.
By the quadratic expansion (Proposition~1),
$\KL(\widehat p\|q_r)=\KL(p_r^*\|q_r)+\frac{1}{2n}u^\top H_r^* u
+R_3(u/\sqrt n)$ with
$|R_3(u/\sqrt n)|\le 2\|u\|_2^3/\{3(p_{\min,r}^*)^2 n^{3/2}\}$, the same
constant $2/\{3(p_{\min,r}^*)^2\}$ as in Proposition~1.  Inside
$\|u\|_2\le\rho_n$ the cubic term contributes
$\exp\{-nR_3\}=1+O_r(\rho_n^3/\sqrt n)$, so the good-type weight is
$\exp\{-\tfrac12 u^\top H_r^* u\}\{1+O_r(\rho_n^3/\sqrt n+n^{-1})\}$, the
relative error and its constant depending on the fixed chart.\qed

\section{Proof of Theorem~2 (Gaussian Localization on the Constraint Chart)}%
\label{sec:proof-gdf}

Decompose the constrained predictive law as
\[
\mu_{n,m}^{(r)}(A)=S_n^{(\mathrm{in})}(A)+S_n^{(\mathrm{out})}(A),
\]
where $S_n^{(\mathrm{in})}$ sums over good types
$\widehat{p}\in\mathcal{G}_{n,r}$ and $S_n^{(\mathrm{out})}$ over bad
types $\widehat{p}\in\mathcal{B}_{n,r}$.

By Lemma~1,
$|S_n^{(\mathrm{out})}(A)|\le O\{n^{-1}(\log n)^{s_r/2}\}$.

For the inner sum, Lemma~2 replaces each type weight by its Gaussian
approximation.  Two operations remain: converting the lattice sum into an
integral and extending that integral beyond the localization ball.

\medskip\noindent\textit{Lattice-to-integral conversion (coarea).}
The realized types lie on the simplex lattice
$\Pcal_n=\{p:np\in\mathbb{Z}_{\ge0}^{k_r},\ \sum_j p_j=1\}$; they do
\emph{not} form the standard lattice $n^{-1/2}\mathbb{Z}^{k_r-1-d}$ in an
arbitrary orthonormal tangent basis.  We therefore work in the affine
simplex chart and treat the moment condition as the constraint surface
$\{A_r p=\alpha\}$.  Under lattice regularity the exactly feasible types,
written in the scaled tangent coordinate
$u=\sqrt{n}\,V_r^\top(\widehat p-p_r^*)$, form a translate of a rank-$s_r$
lattice $\Lambda_r$ of mesh $h=O(n^{-1/2})$ and covolume
  $\mathrm{covol}(\Lambda_r)=c_r\,n^{-s_r/2}$, the reciprocal lattice density
of $\Pcal_n$ restricted to the constraint surface.  Lemma~A4 (\S\,A8) then
converts the constrained lattice sum into an integral,
\[
\mathrm{covol}(\Lambda_r)\!\!\sum_{u\in\Lambda_r+a}\!\! f(u)
=\int_{\R^{s_r}}f(u)\,du+O_r\{h\,M_h(f)\},
\qquad h=n^{-1/2},
\]
where $M_h(f)$ is the summable cellwise derivative envelope of Lemma~A4.
The lemma is applied to the Gaussian kernel and, after a smooth cutoff
outside the localization ball, to the Gaussian-weighted predictive and
Stirling factors.  Their derivative envelopes are $O_r(1+m)$ uniformly in
$n$.  The covolume factor $c_r
n^{-s_r/2}$ is common to $\mathcal{G}_{n,r}$ and to its Gaussian integral,
and therefore \emph{cancels} in the normalized weights $w_{n,r}$; the
residual Riemann-sum error is $O_r\{(1+m)n^{-1/2}\}$ for fixed $r$, and the
truncation at $\|u\|_2\le\rho_n$ is exponentially small.  The constant is
$r$-dependent; no dimension-free bound is claimed.

\medskip\noindent\textit{Tail extension.}
Extending the Gaussian integral from $\|u\|_2\le\rho_n$ to all of
$\R^{s_r}$ adds mass
$\Pbb_{\Nor(0,(H_r^*)^{-1})}(\|u\|_2>\rho_n)
\le C_r(1+\rho_n)^{s_r}
\exp\{-\lambda_{\min}(H_r^*)\rho_n^2/2\}
=O_r\{n^{-2}(\log n)^{s_r/2}\}$.

\medskip\noindent\textit{Exact feasibility and the thickening extension.}
In the exact-feasibility case ($\eta_n=0$, under lattice regularity along
the cofinal set $n\in\mathcal{N}_r$) every feasible type lies on
$\{A_r p=\alpha\}=p_r^*+T_r^*$; the continuous constrained minimizer is
$p_r^*$, generally not a lattice point, and the discrete minimizer over
exactly feasible types, $\bar p_{n,r}$, lies within $O_r(n^{-1})$ of it, so
all remaining arguments proceed on the constraint surface with no normal
displacement.
For the generic thickened case ($\eta_n>0$) let the \emph{continuous} slab
minimizer be
$p_r^{\eta_n}=\argmin\{\KL(p\|q_r):p\in\Delta^{k_r-1},\ \|A_r p-\alpha\|
\le\eta_n\}$; under the relative-interior and full-rank conditions of
Definition~1 the Karush--Kuhn--Tucker conditions and the implicit-function
theorem give $\|p_r^{\eta_n}-p_r^*\|=O_r(\eta_n)$, the moment-normal
displacement forced by the tolerance, while $n\eta_n\to\infty$ makes
$E_{n,r}$ nonempty.  The feasible types then carry a normal component of
size $O(\eta_n)$, so the \emph{sharp} reduced-rank Gaussian of Theorem~2
would require a two-scale conditional local-limit theorem (tangent scale
$n^{-1/2}$, active normal scale $n^{-1}$), which we do not develop; the
collapse rate in this case is supplied instead by Theorem~4, proved in
\S\,A5.1 without any tangent expansion.

Combining, in the exact-feasibility case, the local Taylor/Stirling error
$O_r\{(1+m)\rho_n^3/\sqrt{n}\}$, the quadrature error
$O_r\{(1+m)n^{-1/2}\}$, the collision
error $O(m^2/n)$ (\S\,A5) and the exponentially small Gaussian tail gives
the remainder $R_{n,m,r}$ with
$\sup_A|R_{n,m,r}(A)|\le C_r\{(1+m)n^{-1/2}(\log n)^{3/2}+m^2/n\}$,
the constant depending on the fixed chart; the
$n^{-1/2}(\log n)^{3/2}$ term is not absorbed into $o(n^{-1/2})$.  The
thickened extension adds $O_r(m\eta_n)$.\qed

\section{Proof of Theorem~3 (Exchangeable-Block Collapse)}%
\label{sec:proof-collapse}

The total-variation distance decomposes as
\begin{align*}
\|\mu_{n,m}^{(r)}-(p_r^*)^{\otimes m}\|_{\TV}
&\le\int_{\|v\|_2\le\rho_n}
\|(p_r^*(v))^{\otimes m}-(p_r^*)^{\otimes m}\|_{\TV}
  \;\varphi_{H_r^*}(v)\,dv+|\mathcal{E}_n|\\
&\le m\int_{\|v\|_2\le\rho_n}\|p_r^*(v)-p_r^*\|_{\TV}
  \;\varphi_{H_r^*}(v)\,dv+|\mathcal{E}_n|,
\end{align*}
using the product coupling inequality
$\|P^{\otimes m}-Q^{\otimes m}\|_{\TV}\le m\|P-Q\|_{\TV}$.

The truncated integral is bounded by the corresponding full Gaussian
moment,
\[
m\cdot\Ebb_\varphi\{\|V_r v/\sqrt{n}\|_1\}
\le \frac{m\sqrt{k_r}}{\sqrt{n}}\,\Ebb_\varphi\{\|v\|_2\}
\le \frac{m}{\sqrt{n}}\bigl\{k_r\,\tr(H_r^{*-1})\bigr\}^{1/2},
\]
using $\|V_r v\|_1\le\sqrt{k_r}\,\|v\|_2$ (orthonormal $V_r$) and
$\Ebb_\varphi\{\|v\|_2\}\le\{\tr(H_r^{*-1})\}^{1/2}$.  This is the leading
term of Theorem~3; it carries \emph{no} logarithmic factor, being an
expectation of $\|v\|_2$ under the fixed Gaussian.

The finite-population term comes from the collision bound~(11) of the main
paper: drawing an $m$-block without replacement from the $n$-ensemble,
$\|\mathrm{HG}_{n,p,m}-p^{\otimes m}\|_{\TV}\le 1-(n)_m/n^m\le m(m-1)/(2n)$
uniformly in $p$, giving the term $m(m-1)/n$.  This bound does \emph{not}
involve $1/p_{\min,r}^*$.

The remaining Theorem-2 error, $\sup_A|R_{n,m,r}(A)|\le
C_r(1+m)n^{-1/2}(\log n)^{3/2}$, enters explicitly and is
\emph{not} absorbed into $o(n^{-1/2})$.  Collecting yields the bound of
Theorem~3 in the exact-feasibility case; the thickened extension adds
$O_r(m\eta_n)$.\qed

\subsection{Window-Robust Collapse (Theorem~4)}\label{sec:proof-window}

For generic charts we drop lattice regularity entirely.  The argument uses
convexity and type counting only: no tangent expansion, no Gaussian
localization and no covolume cancellation enter, so neither Proposition~1
nor Lemmas~1--2 is invoked.

\medskip\noindent\textit{Step 1: a Bregman--Pinsker lower bound on the slab.}
Let $\mathcal{S}_{n,r}=\{p\in\Delta^{k_r-1}:\|A_r p-\alpha\|\le\eta_n\}$ be
the continuous slab and $p_r^{\eta_n}$ its $\KL(\cdot\|q_r)$-minimizer
(\S\,A4).  As $\mathcal{S}_{n,r}$ is convex, first-order optimality gives
$\nabla\KL(p_r^{\eta_n}\|q_r)^\top(p-p_r^{\eta_n})\ge0$ for every
$p\in\mathcal{S}_{n,r}$.  The Bregman identity for $\KL$, whose Bregman
divergence in the first argument is again $\KL$, then gives, for all
$p\in\mathcal{S}_{n,r}$,
\begin{equation}\label{eq:bregman}
\begin{split}
\KL(p\|q_r)-\KL(p_r^{\eta_n}\|q_r)
&=\KL(p\|p_r^{\eta_n})
+\nabla\KL(p_r^{\eta_n}\|q_r)^\top(p-p_r^{\eta_n})\\
&\ \ge\ \KL(p\|p_r^{\eta_n})
\ \ge\ \tfrac12\|p-p_r^{\eta_n}\|_1^2 ,
\end{split}
\end{equation}
the last step by Pinsker's inequality.  This holds globally on the slab; it
uses no local expansion and no lattice structure.

\medskip\noindent\textit{Step 2: inward rounding to a feasible type.}
Since $A_rp_r^*=\alpha$, the point $p_r^*$ lies in $\mathcal{S}_{n,r}$ with
zero moment slack, while $\|p_r^{\eta_n}-p_r^*\|_1=O_r(\eta_n)$ (\S\,A4).
For $s=c_r/(n\eta_n)$, which tends to $0$ because $n\eta_n\to\infty$, the
contraction $\tilde p=(1-s)p_r^{\eta_n}+s\,p_r^*$ satisfies
$\|A_r\tilde p-\alpha\|\le(1-s)\eta_n=\eta_n-c_r/n$ and
$\|\tilde p-p_r^{\eta_n}\|_1=s\,O_r(\eta_n)=O_r(n^{-1})$.  Rounding
$n\tilde p$ to integers perturbs the moment slack by $O_r(n^{-1})$, so for
$c_r$ large enough the resulting type $\bar p^{\,\eta}_{n,r}$ lies in
$E_{n,r}$ and obeys $\|\bar p^{\,\eta}_{n,r}-p_r^{\eta_n}\|_1=O_r(n^{-1})$.
This is the only place where $n\eta_n\to\infty$ is used.  As
$\KL(\cdot\|q_r)$ is Lipschitz near $p_r^*$, it follows that
$n\KL(\bar p^{\,\eta}_{n,r}\|q_r)\le n\KL(p_r^{\eta_n}\|q_r)+O_r(1)$, so by
the method of types
$\Pbb(\whatP^{(r)}\in E_{n,r})\ge
(n+1)^{-k_r}\exp\{-n\KL(p_r^{\eta_n}\|q_r)-O_r(1)\}$.

\medskip\noindent\textit{Step 3: concentration.}
There are at most $(n+1)^{k_r}$ types, each of probability at most
$\exp\{-n\KL(p\|q_r)\}$.  Combining this with~\eqref{eq:bregman} and the
denominator bound of Step~2, for every $t>0$,
\[
\Pbb\bigl(\|\whatP^{(r)}-p_r^{\eta_n}\|_1>t \bigm|
\whatP^{(r)}\in E_{n,r}\bigr)
\ \le\ C_r\,(n+1)^{2k_r}\,e^{-nt^2/2}.
\]
Taking $t_n=\{2(2k_r+1)\log(n+1)/n\}^{1/2}=O_r(\sqrt{\log n/n})$ makes the
right-hand side $O(n^{-1})$; since $\|\cdot\|_1\le2$ on the simplex,
integrating the tail gives
\[
\Ebb\bigl\{\|\whatP^{(r)}-p_r^{\eta_n}\|_1 \bigm|
\whatP^{(r)}\in E_{n,r}\bigr\}
=O_r\!\left(\sqrt{\tfrac{\log n}{n}}\right).
\]

\medskip\noindent\textit{Step 4: assembly.}
With $\|p_r^{\eta_n}-p_r^*\|_1=O_r(\eta_n)$, the triangle inequality gives
$\Ebb\{\|\whatP^{(r)}-p_r^*\|_1\mid\whatP^{(r)}\in E_{n,r}\}
=O_r(\eta_n+\sqrt{\log n/n})$.  By the exact identity~(10),
$\mu_{n,m}^{(r)}=\sum_p\mathrm{HG}_{n,p,m}\,w_{n,r}(p)$; the collision
bound~(11) replaces $\mathrm{HG}_{n,p,m}$ by $p^{\otimes m}$ at cost
$m(m-1)/(2n)$ uniformly, hence $m^2/n$ after averaging, and product
coupling gives $\|p^{\otimes m}-(p_r^*)^{\otimes m}\|_{\TV}
\le m\|p-p_r^*\|_{\TV}=\tfrac m2\|p-p_r^*\|_1$.  Averaging over the
conditional law of $\whatP^{(r)}$ yields
\[
\|\mu_{n,m}^{(r)}-(p_r^*)^{\otimes m}\|_{\TV}
\le C_r\Bigl(m\eta_n+m\sqrt{\tfrac{\log n}{n}}+\tfrac{m^2}{n}\Bigr),
\]
which is Theorem~4.  The logarithm is the price of Pinsker-type
concentration in place of the reduced-rank Gaussian of Theorem~3.\qed

\section{Proofs of Parametric and Robustness Results}\label{sec:proof-parametric}

\subsection{Induced Conditional Likelihood and Generalized Posterior}

Under the uniform conditions of Proposition~2, using Theorem~3 in the
exact-lattice case and Theorem~4 for a generic shrinking window, the
constrained block law $\mu_{n,m}^{(r)}(\cdot;\theta)$ collapses to
$(p_{\theta,r}^*)^{\otimes m}$ uniformly over compact $K\subset\Theta$.
The associated product scoring rule is the induced conditional likelihood
$L_r^{\text{\textsc{icl}}}(\theta;X_{1:n})=\prod_{i=1}^n p_{\theta,r}^*(T_r X_i)$
of the main paper, with cell-level log-likelihood
$\ell_r^{\text{\textsc{icl}}}(\theta)=\sum_{i=1}^n\log p_{\theta,r}^*(B_i)$ used in
the implementation.  Updating a prior $\pi$ by this likelihood and the
feasibility indicator $\ind\{\Ccal_{\theta,r}(\alpha)\neq\varnothing\}$
gives, \emph{by definition}, the induced generalized posterior; this is
an exact update, not an approximation of the exact Bayesian posterior of
the conditioning event.  The feasibility indicator is the finite-chart
analogue of the convex-hull condition of \citet{schennach2005}.  In the
ambient setting the analogous product object is
$\prod_{i=1}^n(dP_\theta^*/dQ)(X_i)$.

\subsection{Chart Efficiency of the Induced Conditional-Likelihood
estimator}\label{sec:icl-efficiency}

Because the cell counts are multinomial, criterion~(21) is the exact
log-likelihood for the smooth positive simplex submodel
$\{p_{\theta,r}^*:\theta\in\Theta\}$.  The induced conditional-likelihood estimator is a
maximum-likelihood estimator on this chart; the consistency condition below
excludes distant aliases of the same cell law.  The scope of the resulting
efficiency statement is deliberately narrow.

\begin{proposition}[Chart efficiency]\label{prop:hajek}
Fix a moment-adequate partition $\Pi_r$ and let
Assumption~1 hold.  Suppose the chart model is
correctly specified, that is, the true cell law equals
$p_{\theta_0,r}^*$ for some $\theta_0$ interior to $\Theta$, and that
$J_{\theta_0,r}=\partial p_{\theta,r}^*/\partial\theta^\top$ has full
column rank.  Assume $\theta\mapsto p_{\theta,r}^*$ is twice continuously
differentiable near $\theta_0$, and let
$\widehat\theta_{\text{\textsc{icl}}}$ denote a consistent local maximizer of the
chart likelihood.  Then, as $n\to\infty$ with $r$ fixed:
\begin{enumerate}
\item[(i)] the chart experiment
    $\Ecal_{n,r}=\{\mathrm{Mult}(n,p_{\theta,r}^*):\theta\in\Theta\}$ is
locally asymptotically normal at $\theta_0$ with information
$I_{\theta_0,r}$ of~(18);
\item[(ii)] $\sqrt n(\widehat\theta_{\text{\textsc{icl}}}-\theta_0)
\Rightarrow\Nor(0,I_{\theta_0,r}^{-1})$;
\item[(iii)] $\widehat\theta_{\text{\textsc{icl}}}$ is asymptotically efficient
\emph{within the chart experiment}: by the convolution theorem of
\citet{hajek1970}, the limit law of any regular estimator sequence in
$\Ecal_{n,r}$ is $\Nor(0,I_{\theta_0,r}^{-1})*M$ for some distribution
$M$, and the local asymptotic minimax bound is attained.
\end{enumerate}
\end{proposition}

\begin{proof}
By Definition~1(iii) the projected cell law is interior, so
$p_{\theta,r}^*\ge p_{\min,r}^*>0$ on a neighbourhood of $\theta_0$.  By
the exponential-family duality of Lemma~A1,
$p_{\theta,r}^*\propto q_r\exp\{\lambda(\theta)^\top h_r\}$ with
$\lambda(\theta)$ continuously differentiable, the implicit function
theorem applying because $\widetilde A_r$ has full row rank and the dual
moment covariance $\Omega_{\theta,r}$ is nonsingular.  Hence
$\theta\mapsto p_{\theta,r}^*$ is $C^1$ with values bounded away from
zero, so $\theta\mapsto\sqrt{p_{\theta,r}^*}$ is $C^1$ and the submodel is
differentiable in quadratic mean, with score
$s_\theta(j)=\partial_\theta\log p_{\theta,r,j}^*
=(\partial_\theta\lambda)^\top\{h_r(j)-\Ebb_{p_{\theta,r}^*}h_r\}$
and information $I_{\theta,r}=J_{\theta,r}^\top\diag(1/p_{\theta,r}^*)
J_{\theta,r}$, nonsingular by the rank hypothesis.  Local asymptotic
normality of the multinomial submodel and (ii)--(iii) are then the
standard consequences \citep[Thms.~7.2, 8.8 and 8.11]{vanderVaart1998}; no
property of the block law, and no finite-population correction, is used.
\end{proof}

\medskip\noindent\textit{Scope of the efficiency statement.}  Two restrictions are
essential, and neither is cosmetic.

\emph{Efficiency is within the chart.}  The optimality in
Proposition~\ref{prop:hajek}(iii) is relative to estimators that observe
only the cell counts, that is, within $\Ecal_{n,r}$.  It is not a claim
of efficiency in the ambient experiment: the chart discards information,
and \S\,\ref{sec:lecam-supp} is precisely the quantification of how much.
At sample size $n$, $\delta(\Ecal_n,\Ecal_{n,r})=0$: an ambient observation
can reproduce a chart procedure exactly by applying the partition map.
The reverse deficiency $\delta(\Ecal_{n,r},\Ecal_n)=O(n\omega_r)$ controls,
for bounded loss, how well a chart observation can approximate an ambient
procedure.  We assert no quantitative efficiency-loss bound for quadratic
risk.

\emph{Correct specification on the chart is a hypothesis, not a
consequence.}  The true cell law is the partition image of the ambient
projected law, $T_{r\#}P_{\theta_0}^*$, whereas the model member is the
projection of the partitioned baseline, $p_{\theta_0,r}^*$; projection and
partitioning do not commute.  Lemma~A3 bounds the one-block discrepancy by
$O(\omega_r)$, and Corollary~1 lifts the predictive comparison to the
original sample space with product error $O(m\omega_r)$.  The two relevant
regimes are therefore
\[
\begin{array}{lll}
\text{inference centred at }\theta_0
  &:\quad \sqrt n\,\omega_{r_n}\to0
  &\text{(an undersmoothed chart)},\\
\text{predictive calibration}
  &:\quad k_{r_n}\asymp n^{1/3},\quad \omega_{r_n}\asymp n^{-1/3}
  &\text{(the balance in Remark~6)}.
\end{array}
\]
In the second regime $\sqrt n\,\omega_{r_n}\asymp n^{1/6}\to\infty$, so
the chart discrepancy dominates the sampling error.
Thus Proposition~\ref{prop:hajek} is fixed-$r$, as are Theorems~2--3
(Remark~5).  Under refinement the centring moves to the chart pseudo-true
value, and the relevant second-order object is the misspecification curvature $K_r$ of
\S\,\ref{sec:proof-parametric}, not $I_{\theta,r}$.

\subsection{Proof of Proposition~3 (Projection-Induced Moment Criterion)}

Fix the cell-moment matrix $A_r$ and let the target $\alpha(\theta)$ be
twice continuously differentiable.  Under the chart-correct law
$p_{\theta_0,r}^*$, the multinomial central limit theorem gives
$A_r\widehat p_n-\alpha(\theta_0)=O_\Pbb(n^{-1/2})$.  By the assumed
relative-interior condition, with probability tending to one there is a
finite $\widehat\lambda_r$ satisfying
$\nabla\psi_r(\widehat\lambda_r)=A_r\widehat p_n$.
Since $p_{\theta,r,j}^*=q_{r,j}\exp\{\lambda_{\theta,r}^\top
h_r(j)-\psi_r(\lambda_{\theta,r})\}$ with log-partition
$\psi_r(\lambda)=\log\sum_j q_{r,j}\exp\{\lambda^\top h_r(j)\}$, the induced
log-likelihood is \emph{exactly}
\[
-\tfrac1n\ell_r^{\text{\textsc{icl}}}(\theta)
=-\sum_j\widehat p_{n,j}\log q_{r,j}
-\lambda_{\theta,r}^\top\widehat\alpha_n+\psi_r(\lambda_{\theta,r}),
\qquad \widehat\alpha_n=A_r\widehat p_n,
\]
whose first term is free of $\theta$.  The tilt is the exponential-family
dual variable $\nabla\psi_r(\lambda_{\theta,r})=\alpha(\theta)$, and
$\nabla^2\psi_r(\lambda)=A_r\Sigma_{\lambda,r}A_r^\top$ is the cell-moment
covariance under the tilt, equal at $\theta$ to
$\Omega_{\theta,r}=A_r\Sigma_{\theta,r}A_r^\top$,
$\Sigma_{\theta,r}=\diag(p_{\theta,r}^*)-p_{\theta,r}^*p_{\theta,r}^{*\top}$.
Expanding the smooth
convex dual loss
$\psi_r(\lambda_{\theta,r})-\lambda_{\theta,r}^\top\widehat\alpha_n$ about
$\widehat\lambda_r$, with $\lambda_{\theta,r}-\widehat\lambda_r
=\Omega_{\theta,r}^{-1}\{\alpha(\theta)-\widehat\alpha_n\}
+O(\|\alpha(\theta)-\widehat\alpha_n\|^2)$, gives uniformly over
$\|\theta-\theta_0\|_2\le Cn^{-1/2}$
\[
-\log L_r^{\text{\textsc{icl}}}(\theta)
=(\theta\text{-free})
+\frac n2\{\widehat\alpha_n-\alpha(\theta)\}^\top\Omega_{\theta,r}^{-1}
\{\widehat\alpha_n-\alpha(\theta)\}+O_\Pbb(n^{-1/2}).
\]
This is the stated quadratic criterion with weight $\Omega_{\theta,r}^{-1}$,
obtained directly by Legendre duality, without decomposing $\delta_\theta$
into pinned and free tangent parts.  Differentiating
$\nabla\psi_r(\lambda_{\theta,r})=\alpha(\theta)$ gives the parameter
information
$I_{\theta,r}=\dot\alpha(\theta)^\top\Omega_{\theta,r}^{-1}\dot\alpha(\theta)$.\qed

On the empirical-support chart, with $p_{\theta,r}^*$ an exponential tilt
of $\whatP$, $\Omega_{\theta,r}$ is the covariance of the fixed cell
moments under the tilted empirical law; under covariance regularity,
differentiability, local identification and stochastic equicontinuity it
converges to the long-run moment covariance and the criterion matches its
efficient generalized method-of-moments counterpart for the fixed moment
map.

For a general parameter-dependent moment map $g(x,\theta)$ the cell-moment
matrix $A_{\theta,r}$ itself varies with $\theta$; differentiating
$A_{\theta,r}p^*_{\theta,r}=\alpha(\theta)$ then produces an additional
$\dot A_{\theta,r}p^*_{\theta,r}$ term and a moving tangent space.  A
treatment would require uniform-in-$\theta$ control of the cell-moment
matrices and their derivatives, beyond the fixed-matrix expansion here.

\subsection{Misspecification: Curvature and Sandwich Covariance}

The demoted misspecification paragraph of the main paper rests on the
following standard M-estimation and generalized-Bayes calculation, which
we record but do not develop into a theorem.  Let
$R_r(\theta)=-\Ebb_{P_0}\log p_{\theta,r}^*(T_r X)$, with minimizer
$\theta_{0,r}$ and parameter curvature
$K_r=\nabla_\theta^2 R_r(\theta_{0,r})$, assumed positive definite.  A
second-order expansion of the empirical risk about $\theta_{0,r}$ gives a
quadratic log-posterior with Hessian $nK_r+o_\Pbb(n)$, so that under the
usual regularity conditions
\[
\Pbb(\|\theta-\theta_{0,r}\|>\varepsilon\mid X_{1:n})
\le C\exp\{-c\,n\,\lambda_{\min}(K_r)\varepsilon^2\},
\]
governed by $\lambda_{\min}(K_r)$, not $\lambda_{\min}(H_{\theta_0,r}^*)$.
The associated estimator has sandwich covariance
$n^{-1}K_r^{-1}\Sigma_r K_r^{-1}$ with
$\Sigma_r=\mathrm{Var}_{P_0}\{\nabla_\theta\log p_{\theta_{0,r},r}^*(T_r X)\}$.
For each fixed $\theta$, the predictive counterpart follows from
Theorem~3, or Theorem~4 for a generic window, with limiting law
$(p_{\theta,r}^*)^{\otimes m}$.  Uniform collapse near $\theta_{0,r}$
together with the local Lipschitz bound
$\|p_{\theta,r}^*-p_{\theta_{0,r},r}^*\|_{\TV}
\le L_r\|\theta-\theta_{0,r}\|$
then calibrates plug-in prediction.  Thus parameter concentration is
governed by $K_r$, whereas predictive concentration on the constraint
chart is governed by $H_r^*$.  The broader framework belongs to the
generalized-posterior literature and is not re-derived here.

\section{The Le~Cam Bridge}\label{sec:lecam-supp}

The finite-partition model is a coordinate realization of the
measure-theoretic construction; the deficiency calculation below
quantifies the relationship.

\subsection{Deficiency and the Discretization Gap}

Let $\Ecal_n=\{P_\theta^{\otimes n}:\theta\in\Theta\}$ denote the
original experiment on $\Xcal^n$ and let
$T_r:\Xcal\to\{1,\ldots,k_r\}$ be the partition map.  The discretized
experiment is
$\Ecal_{n,r}=\{(P_\theta\circ T_r^{-1})^{\otimes n}:\theta\in\Theta\}$.

Since $T_r$ is deterministic, the original experiment can always simulate
the discretized one.  The partition map $T_r^{\otimes n}$ is the required
Markov kernel, and hence $\delta(\Ecal_n,\Ecal_{n,r})=0$.  The reverse
direction contains the substance of the comparison.

\subsection{Reverse Kernel and Approximation Rate}

The nontrivial direction is $\delta(\Ecal_{n,r},\Ecal_n)$: how well
the discretized experiment can reconstruct the original.  Introduce the
reverse kernel
$\kappa_{r,j}(dx)=Q(dx\mid A_j^{(r)})$, mapping each cell to the
baseline conditional distribution within it.  This produces the lifted
law
$\widetilde{P}_{\theta,r}(dx)
=\sum_{j=1}^{k_r}P_\theta(A_j^{(r)})\,\kappa_{r,j}(dx)$.

\begin{assumption}[Partition regularity]\label{ass:partition}
The sample space $(\Xcal,\rho)$ is a metric space and the refining
sequence $\{\Pi_r\}_{r\ge 1}$ satisfies:
\begin{enumerate}
\item $\Pi_{r+1}$ refines $\Pi_r$;
\item $\sigma(\bigcup_r\Pi_r)=\Acal$;
\item the mesh diameter
  $\mathrm{mesh}_r=\sup_j\mathrm{diam}_\rho(A_j^{(r)})\to 0$ as $r\to\infty$,
  the diameter taken in the metric $\rho$.  This is distinct from the cell
  oscillation $\omega_r$ of Assumption~2; for Lipschitz $h$ the two are
  related by $\omega_r\le L_h\,\mathrm{mesh}_r$.
\end{enumerate}
\end{assumption}

On an unbounded sample space a finite partition cannot have globally
vanishing mesh, since the tail cells are unbounded; condition~(iii) is then
understood on a compact domain, or after truncation, with the unbounded
tail cells controlled directly by the uniform cell-regularity hypothesis
below rather than by the mesh.

\begin{theorem}[Le~Cam approximation under refinement]\label{thm:lecam}
Let Assumptions~1 and~2 of the main paper and
Assumption~\ref{ass:partition} hold, with $\omega_r$ the cell oscillation of
Assumption~2.  Suppose that $h$ is Lipschitz with constant $L_h$ and that $P_\theta$ has density
$f_\theta=dP_\theta/dQ$ satisfying
$\|f_\theta-f_{\theta'}\|_{L^1(Q)}\le C_f\|\theta-\theta'\|$ for
$\theta,\theta'$ in a compact parameter set $\Theta_0$.  Assume
further, as a direct condition, that the densities are uniformly regular
on cells,
\[
\sup_{\theta\in\Theta_0}\,
\sup_{x,x'\in A_j^{(r)}}|f_\theta(x)-f_\theta(x')|
\le C_f'\,\omega_r\quad\text{for all }j.
\]
Then the Le~Cam
distance satisfies
\[
\Delta(\Ecal_n,\Ecal_{n,r})\le n\cdot
\sup_{\theta\in\Theta_0}\TV(P_\theta,\widetilde{P}_{\theta,r})
=O(n\omega_r).
\]
In particular,
$\Delta(\Ecal_n,\Ecal_{n,r})\to 0$ as $r\to\infty$ for fixed $n$.
\end{theorem}

On a bounded domain the cell-regularity condition follows from Lipschitz
$h$ and a bounded tilt $\lambda_\theta$.  On an unbounded domain it is a
genuine additional assumption: a small cell oscillation of $h$ does not by
itself control the tail oscillation of $f_\theta$.

\begin{proof}
The Le~Cam distance satisfies
$\Delta(\Ecal_n,\Ecal_{n,r})\le n\cdot\delta_1$ by tensorization of
Markov kernels, where $\delta_1$ is the single-observation deficiency.
We have $\delta(\Ecal_1,\Ecal_{1,r})=0$ (the original simulates the
discretized via $T_r$).  For the reverse direction, the reverse kernel
$\kappa_r$ maps the discretized observation $j$ to
$\kappa_{r,j}(\cdot)=Q(\cdot\mid A_j^{(r)})$.

The total variation between $P_\theta$ and its reconstruction
$\widetilde{P}_{\theta,r}$ is
\[
\TV(P_\theta,\widetilde{P}_{\theta,r})
=\tfrac{1}{2}\int_\Xcal
\left|f_\theta(x)-\sum_{j=1}^{k_r}
\frac{P_\theta(A_j^{(r)})}{Q(A_j^{(r)})}
\ind_{A_j^{(r)}}(x)\right|Q(dx).
\]
On each cell $A_j^{(r)}$ the integrand is
$|f_\theta(x)-\Ebb_Q\{f_\theta\mid\Pi_r\}(x)|$, a martingale
increment.  The conditional expectation
$\Ebb_Q\{f_\theta\mid\Pi_r\}\to f_\theta$ as $r\to\infty$, both
$Q$-a.e.\ (by the martingale convergence theorem) and in $L^1(Q)$
(by L\'evy's upward theorem).  Here $f_\theta\in L^1(Q)$ because it is a
probability density; the conditional-expectation martingale is therefore
uniformly integrable.

For the rate, the uniform cell-regularity condition bounds the integrand
by $C_f'\,\omega_r$ on every cell; integrating over all cells yields
$\TV(P_\theta,\widetilde{P}_{\theta,r})=O(\omega_r)$ uniformly over
$\theta\in\Theta_0$, and tensorization gives
$\Delta(\Ecal_n,\Ecal_{n,r})=O(n\omega_r)$.
\end{proof}

For prediction we do \emph{not} claim a direct conditional-deficiency
bound between $\mu_{n,m}$, the ambient conditioned block law, and its chart
counterpart $\mu_{n,m}^{(r)}$: conditioning can amplify an unconditioned
discrepancy.  The loop is closed through the projected
laws instead.  By the projection-approximation lemma (\S\,A8,
Lemma~A3), $\TV(\widetilde P_r^*,P^*)=O(\omega_r)$; lifting the chart law
through the reverse kernel, the block limits $(P^*)^{\otimes m}$ and
$(\widetilde P_r^*)^{\otimes m}$, both on $\Xcal^m$, differ by
$O(m\omega_r)$, while the cell probabilities of the two projected laws
differ by $O(\omega_r)$ under the same lemma.  Applying the reverse kernel
to the entire chart predictive law and using total-variation contraction
closes the ambient/finite loop in Corollary~1 of the main paper.

\section{Auxiliary Results}\label{sec:auxiliary}

\begin{lemma}[Existence and uniqueness of the I-projection]%
\label{lem:iprojection}
If $q_{r,j}>0$ for all $j$ and
$\alpha\in\mathrm{ri}\,\mathrm{conv}\{h_r(1),\ldots,h_r(k_r)\}$
(Definition~1), then the finite-dimensional projection $p_r^*$ exists, is
unique, and satisfies $p_{r,j}^*>0$ for all $j$.
\end{lemma}

\begin{proof}
Existence, uniqueness and strict positivity follow from convex duality
for the finite exponential family, not from compactness of the level sets
of $\KL(\cdot\|q_r)$ in the open simplex
$(\Delta^{k_r-1})^\circ=\{p:p_j>0,\;\sum p_j=1\}$, on which those level
sets need not be compact: a minimizing sequence may approach the boundary
while $\KL$ stays finite.  On the simplex, the primal problem
$\min\{\KL(p\|q_r):A_rp=\alpha\}$ is strictly convex in $p$.  Normalization
is absorbed by the log-partition function, so its identifiable dual is the
$d$-dimensional smooth concave program
$\max_{\lambda\in\R^d}\{\lambda^\top\alpha-\psi_r(\lambda)\}$, where
$\psi_r(\lambda)=\log\sum_j q_{r,j}\exp\{\lambda^\top h_r(j)\}$.
The relative-interior condition
$\alpha\in\mathrm{ri}\,\mathrm{conv}\{h_r(1),\ldots,h_r(k_r)\}$ places the
target in the interior of the mean-value domain, so by steepness of
$\psi_r$ (\S\,A1) the dual has a finite maximizer $\lambda_r^*$.  Strong
duality then gives the primal optimum in Gibbs form
$p_{r,j}^*=q_{r,j}\exp\{\lambda_r^{*\top}h_r(j)-\psi_r(\lambda_r^*)\}>0$,
feasible, unique by strict primal convexity, and strictly positive since
$q_{r,j}>0$.  A target that is merely feasible, lying on the boundary of
the mean domain, need not admit a finite dual maximizer and so need not
yield strict positivity.
\end{proof}

\begin{lemma}[Curvature bounds and numerical stability]%
\label{lem:curvature}
The reduced Hessian satisfies:
\begin{enumerate}
\item $1/\tr\{(H_r^*)^{-1}\}\le\lambda_{\min}(H_r^*)\le s_r/\tr\{(H_r^*)^{-1}\}$, where $s_r=k_r-1-d$;
\item $\tr(H_r^*)=\sum_{j=1}^{k_r}\|V_{r,j}\|_2^2/p_r^*(j)$, where
  $V_{r,j}$ denotes the $j$th row of $V_r$;
\item if $\hat{p}_r^*$ is a numerical approximation to $p_r^*$ with
  $\|\hat{p}_r^*-p_r^*\|_\infty\le\varepsilon$, then
  $|\hat{\lambda}_{\min}-\lambda_{\min}(H_r^*)|
  \le 2\varepsilon/(p_{\min,r}^*)^2+O(\varepsilon^2)$.
\end{enumerate}
\end{lemma}

\begin{proof}
For part~(i), $\tr\{(H_r^*)^{-1}\}=\sum_i 1/\lambda_i$; the single term
$1/\lambda_{\min}$ gives $\tr\{(H_r^*)^{-1}\}\ge 1/\lambda_{\min}$, while
bounding every term by $1/\lambda_{\min}$ gives
$\tr\{(H_r^*)^{-1}\}\le s_r/\lambda_{\min}$, and rearranging yields the
two-sided bound.  Part~(ii) is a direct computation from the
definition $H_r^*=V_r^\top\diag(1/p_r^*)V_r$.  For part~(iii),
Weyl's perturbation theorem gives
$|\hat{\lambda}_{\min}-\lambda_{\min}(H_r^*)|
\le\|H_r^*-\hat{H}_r^*\|_2$, and the spectral norm of the
difference $V_r^\top\{\diag(1/p_r^*)-\diag(1/\hat{p}_r^*)\}V_r$ is
bounded by
$\max_j|1/p_{r,j}^*-1/\hat{p}_{r,j}^*|
\le\varepsilon/(p_{\min,r}^*-\varepsilon)^2$,
which gives the stated bound for small $\varepsilon$.
\end{proof}

\begin{lemma}[Projection approximation]\label{lem:spectral}
Let Assumptions~1 and~\ref{ass:partition} hold, with $h$ Lipschitz and
$\mathrm{Cov}_{P^*}\{h\}$ nonsingular.  Suppose $\lambda^*$ has a compact
neighbourhood $K\subset\mathrm{int}\,\mathrm{dom}\,\psi$ on which
\[
\int (1+\|h(x)\|)\,
\exp\{\sup_{\lambda\in K}\lambda^\top h(x)\}\,Q(dx)<\infty .
\]
This local exponential-envelope condition is automatic for bounded $h$,
in particular after compact truncation.  Then the discrete dual moment map
converges uniformly to the ambient
one, so by the implicit-function theorem
$\|\lambda_r^*-\lambda^*\|=O(\omega_r)$, and the lifted projection
$\widetilde P_r^*=\sum_j p_{r,j}^*\,\kappa_{r,j}$ satisfies
$\TV(\widetilde P_r^*,P^*)=O(\omega_r)$.  In particular the cell
probabilities of the ambient projection and the finite I-projection
differ by $|p_{r,j}^*-P^*(A_j^{(r)})|=O(\omega_r)$; they are \emph{not}
assumed equal.
\end{lemma}

\begin{proof}
By the oscillation condition the cell representative $h_r(j)=\Ebb_Q\{h\mid
A_j^{(r)}\}$ differs from $h(x)$ on $A_j^{(r)}$ by $O(\omega_r)$.
The mean-value theorem and the displayed envelope therefore give
$\sup_{\lambda\in K}\{|\psi_r(\lambda)-\psi(\lambda)|
+\|\nabla\psi_r(\lambda)-\nabla\psi(\lambda)\|\}=O(\omega_r)$.
Since $\nabla^2\psi(\lambda^*)=
\mathrm{Cov}_{P^*}\{h\}$ is nonsingular, the implicit-function theorem
applied to $\nabla\psi_r(\lambda_r^*)=\alpha$ gives
$\|\lambda_r^*-\lambda^*\|=O(\omega_r)$.  Comparing the tilted densities
cell by cell and using the same envelope to integrate terms involving
$\|h\|$ gives $\TV(\widetilde P_r^*,P^*)=O(\omega_r)$; the
cell-probability bound follows from $\widetilde P_r^*(A_j^{(r)})=p_{r,j}^*$.
Partition-invariant convergence of the \emph{unnormalized} eigenvalues of
$H_r^*$ is a separate matter and is left open (main paper, \S\,9), since
on a balanced chart $H_r^*$ scales like $k_r I$.
\end{proof}

\begin{lemma}[Lattice Riemann sum]\label{lem:riemann}
Let $\Lambda_h+a\subset\R^{s_r}$ be a translated full-rank lattice whose
fundamental cells $\{C_u:u\in\Lambda_h+a\}$ have volume
$\mathrm{covol}(\Lambda_h)\asymp h^{s_r}$, contain their indexed lattice
points, and have diameter at most $c_rh$.  For $f\in C^1(\R^{s_r})$ define
\[
M_h(f)=\sum_{u\in\Lambda_h+a}\mathrm{vol}(C_u)
       \sup_{v\in C_u}\|\nabla f(v)\|_2 .
\]
If $M_h(f)<\infty$, then
\begin{equation}\label{eq:riemann}
\left|\mathrm{covol}(\Lambda_h)\!\!
\sum_{u\in\Lambda_h+a}\!\! f(u)-\int_{\R^{s_r}}f(u)\,du\right|
\le c_rh\,M_h(f).
\end{equation}
In particular, if
$\|\nabla f(v)\|_2\le C(1+\|v\|_2)^{-s_r-1}$, uniformly over the functions
under consideration, then $\sup_{0<h\le h_0}M_h(f)<\infty$ and the error is
$O_r(h)$.

For the applications in \S\S\,A3--A4, the Gaussian kernel has this decay
directly.  The Gaussian-weighted Stirling and predictive kernels are first
multiplied by a smooth cutoff that equals one on
$\{\|u\|_2\le\rho_n\}$ and vanishes before any coordinate of
$p_r^*+V_ru/\sqrt n$ can reach zero.  Their cellwise derivative envelopes
are $O_r(1+m)$; the cutoff remainder is a Gaussian tail.  Thus, with
$h=n^{-1/2}$, the quadrature error is
$O_r\{(1+m)n^{-1/2}\}$, which is absorbed by the remainder in Theorem~2.
\end{lemma}

\begin{proof}
Since $u\in C_u$ and $\mathrm{diam}(C_u)\le c_rh$, the mean-value theorem
gives
$|f(u)-f(v)|\le c_rh\sup_{w\in C_u}\|\nabla f(w)\|_2$ for $v\in C_u$.
Consequently,
\begin{align*}
\left|\mathrm{covol}(\Lambda_h)\!\!
\sum_{u\in\Lambda_h+a}\!\! f(u)-\int_{\R^{s_r}}f\right|
&\le\sum_u\int_{C_u}|f(u)-f(v)|\,dv\\
&\le c_rh\sum_u\mathrm{vol}(C_u)
       \sup_{w\in C_u}\|\nabla f(w)\|_2
=c_rh\,M_h(f).
\end{align*}
For a fixed lattice shape, summing the polynomial envelope over cells
gives a bounded Riemann sum whenever its decay exponent exceeds $s_r$.
Gaussian times polynomial envelopes satisfy this condition.  On the cutoff
support, all cell probabilities are at least $p_{\min,r}^*/2$ for large
$n$; differentiating an $m$-fold product probability costs at most a
fixed-chart multiple of $m/\sqrt n$, while differentiating the Gaussian
and cutoff costs a fixed-chart polynomial factor.  This verifies the
stated $O_r(1+m)$ envelope and completes the application.
\end{proof}

\section{Simulation Design and Implementation Details}%
\label{sec:sim-details-supp}

Section~8 of the main paper contains two complementary calculations: an
exact enumeration of the conditioned block law on a three-cell chart, and
a parametric Gaussian study of geometry and misspecification.  Both are
illustrative rather than exhaustive.

\subsection{Exact Enumeration of the Conditioned Block Law}

For Fig.~2 of the main paper, take the three-cell baseline
$q=(0{.}2,0{.}5,0{.}3)$, cell moments $h=(-1,0,1)$ and target $\alpha=0$.
The Gibbs multiplier is
$\lambda^*=\tfrac12\log(q_1/q_3)$, giving
$p^*=(0{.}24744871,0{.}50510257,0{.}24744871)$.  Exact feasibility requires
equal counts in cells one and three, so every feasible population has
counts $(t,n-2t,t)$ for $0\le t\le\lfloor n/2\rfloor$.  Its conditional
weight is
\[
w_n(t)\ \propto\
\frac{n!}{t!^2(n-2t)!}\,q_1^t q_2^{\,n-2t}q_3^t .
\]
For an ordered block $z_{1:m}$ with cell counts $c_j(z_{1:m})$, the exact
conditional probability is therefore
\[
\mu_{n,m}(z_{1:m})
=\sum_{t=0}^{\lfloor n/2\rfloor}w_n(t)
\frac{\prod_{j=1}^3\{N_j(t)\}_{c_j(z_{1:m})}}{(n)_m},
\qquad N(t)=(t,n-2t,t).
\]
We enumerate all $3^m$ ordered blocks for $m\in\{1,2,4\}$ and
$n\in\{40,60,100,160,250,400,650,1000,1600,2500\}$.  The plotted
total-variation distances are consequently exact up to floating-point
rounding.  The deterministic script
\texttt{direct\_block\_collapse.py} and its output table
\texttt{direct\_block\_collapse.csv} are included in the simulation
package.

\subsection{Data-Generating Processes}

Two regimes are considered.  Under Model~C (correct specification),
$X_1,\ldots,X_n\stackrel{\mathrm{iid}}{\sim}N(0,1)$ and the moment
conditions $E(X)=0$, $E(X^2)=1$ are exactly satisfied.  Under
Model~M (misspecification), $X_i\sim N(0,4)$, i.e.\ the data
variance is $\sigma^2=4$ ($\sigma=2$) while the constraint still pins the
variance to one.  The fourfold variance ratio produces a substantial
change of measure and makes the differences among the methods visible.
We use sample sizes $n\in\{1000,5000\}$.
Throughout, $q_{M,r}$ denotes the cell law of the $N(0,4)$ measure on
the $Q$-quantile partition; the baseline cell law $q_r$ is uniform.
The extreme cells of the quantile partition are unbounded, so the
refinement regularity of Assumption~2 of the main paper does not apply to
this fixed chart; the $k_r=20$ computation is a fixed-chart calculation.

\subsection{Partition and Baseline}

The baseline measure $Q$ is $N(0,1)$.  The quantile partition of
size $k_r=20$ is formed by dividing the standard normal into
$20$ equiprobable cells: the breaks are
$\Phi^{-1}(j/20)$ for $j=0,\ldots,20$, with the extreme breaks
set to $\pm\infty$.  Cell probabilities under the baseline and
conditional cell moments (mean and raw second moment within each
cell) are computed in closed form from the truncated-normal formulae.
We collect these moments in the $2\times k_r$ matrix $A_r$ from
\S\,2.4 of the main paper.  Its $j$th column, $h_r(j)$, consists of the
conditional mean and raw second moment in cell~$j$.

\subsection{Estimation Procedures}

The comparison starts from the moment system
$g(x,\theta)=\{x-\theta,\,(x-\theta)^2-1\}$, equivalently the cell
moment map $h(x)=(x,x^2)$ with target
$\alpha(\theta)=(\theta,\theta^2+1)$, which pins the mean to $\theta$
and the variance to one.  The two tilted procedures enforce both
conditions.  The sample-mean benchmark uses only the first, correctly
specified location condition and therefore serves a different descriptive
role.

\textit{Generalized method of moments.}
The location condition $\Ebb\{g_1(X,\theta)\}=\Ebb(X-\theta)=0$ is
just identified and yields the method-of-moments estimator
$\hat\theta_{\text{\textsc{gmm}}}=\bar X$.  We adopt this as the
method-of-moments benchmark.  It identifies $\theta$ from the correctly specified location
moment, leaving the variance condition $g_2$ as an over-identifying
restriction that is not used for point estimation.  This
variance-agnostic comparator is specific to the present comparison.  Its
role is descriptive, not an efficiency claim: using $\bar X$ isolates the
effect of the misspecified variance condition on the other two procedures.

\textit{Exponentially tilted empirical likelihood.}
For each $\theta$ the exponential tilt $\lambda(\theta)$ solves the
saddlepoint condition
$\sum_i w_i(\lambda)\,g(x_i,\theta)=0$, equivalently
$\lambda(\theta)=\argmin_\lambda\sum_i
\exp\{\lambda^\top g(x_i,\theta)\}$, with tilted weights
$w_i(\lambda)=\exp\{\lambda^\top g(x_i,\theta)\}/
\sum_j\exp\{\lambda^\top g(x_j,\theta)\}$.  The Bayesian exponentially
tilted empirical likelihood of \citet{schennach2005} has profile
log-likelihood $\ell_{\text{\textsc{betel}}}(\theta)
=\sum_i\log w_i(\lambda(\theta))$, and
$\hat\theta_{\text{\textsc{betel}}}$ maximizes it over $\theta$ on the grid
below.  As a point estimator this is the ETEL maximizer of the BETEL
profile.  The inner saddlepoint is solved by Newton iteration on
$\lambda$, whose Hessian is the weighted sample covariance
$\sum_i w_i(g_i-\bar g)(g_i-\bar g)^\top$, $g_i=g(x_i,\theta)$.  The
tilt acts on the $n$ \emph{empirical atoms}.  If a few extreme observations
receive most of the weight, the Hessian becomes poorly conditioned and
$\lambda(\theta)$ moves towards the convex-hull boundary, where the target
moment can no longer be reproduced.  The realized sample therefore
determines the difficulty of the inner optimization.

\textit{Induced conditional likelihood.}
The induced conditional-likelihood estimator maximizes
$\ell_r^{\text{\textsc{icl}}}(\theta)
=\sum_{i=1}^n\log p_{\theta,r}^*(B_i)
=\sum_{j=1}^{k_r} n_j\log p_{\theta,r,j}^*$, where $B_i$ is the
partition cell containing $X_i$, $n_j$ is the number of observations
in cell $j$, and $p_{\theta,r}^*$ is the
Kullback--Leibler projection of $q_r$ onto the constraint set
$\{p:\sum_j p_j h_r(j)=\alpha(\theta)\}$.  The projection is computed
by Newton iteration on the Lagrange dual, as described in
Section~2.5 of the main paper, with dual Hessian
$\Omega_{\theta,r}=A_r\{\diag(p_{\theta,r}^*)-p_{\theta,r}^*p_{\theta,r}^{*\top}\}A_r^\top$,
the $d\times d$ cell-moment covariance under the projected law.  This
$d\times d$ dual moment covariance weights the induced quadratic criterion of
\S\,A6.3 and is distinct from the $(k_r\!-\!1\!-\!d)$-dimensional reduced
tangent Hessian $H_{\theta,r}^*$ of \S\,2.4; the two coincide only in
scalar special cases.  Here the tilt acts on the $k_r$ \emph{fixed cells}
rather than on the data: the dual Hessian is bounded and independent of
$n$, the projection exists at every grid value satisfying the
relative-interior feasibility condition, and the data enter only through
the cell counts $n_j$.

The two tilted estimators enforce the same mean--variance system but use
different support.  The method-of-moments benchmark instead applies no tilt,
uses only the location condition and returns $\bar X$.  The tilted empirical
likelihood reweights the observed sample,
so its numerical conditioning depends on the realized data and
deteriorates under misspecification.  The induced conditional likelihood
instead tilts a fixed histogram of bounded cell moments; its conditioning
is determined by the partition.  Under Model~C the tilt is nearly trivial
and the three methods nearly coincide.  Under Model~M, concentration of
the empirical-likelihood weights on extreme observations worsens the
Newton solve and increases dispersion, while the induced-likelihood inner
problem remains governed by the fixed partition, and the estimator's
reported standard deviation is close to that of $\bar X$.

\subsection{Computational Settings}

We evaluate the parameter on the regular grid
$\theta\in\{-0.50,-0.4917,\ldots,0.50\}$ of $121$ points.  All
three estimators are evaluated on this grid, $\hat\theta_{\text{\textsc{gmm}}}$
analytically and the other two by grid search followed by a parabolic
sub-grid refinement (a three-point interpolation about the maximizing grid
point).  The resulting estimates are therefore continuous rather than
grid-valued.  We use $2000$ Monte Carlo replications for the main
comparison (Table~1 of the main paper) and $500$ per sample-size
point for the empirical cell-discrepancy figure.  The convex-hull condition for the tilted
empirical likelihood never failed across the $8000$ replications of the main
comparison.  All runs are seeded for exact reproducibility.

\subsection{Curvature-Refinement Sequence}

The condition-number figure (left panel of Fig.~4 of the main
paper) evaluates
the condition number of the reduced Hessian at the projection for
partition sizes $k_r\in\{4,6,8,10,15,20,30,40,60,80\}$.  Under Model~C the
projection coincides with the baseline (since $Q$ already satisfies the
constraints), and the reduced Hessian $H_r^*$ is $k_r$ times a scalar
matrix, giving $\kappa=1$ at every resolution.  Under Model~M the Model-M
projection $q_{M,r}\to p_{M,r}^*$ differs from the baseline and its reduced
Hessian $H_{M,r}^*$ grows in conditioning, reaching
$\kappa(H_{M,r}^*)\approx 1.379$ at $k_r=80$.

\subsection{Total-Variation Diagnostic}

The mean total-variation distance reported in Table~1 is
$\TV(\hat p_n,p_r^*)=(1/2)\sum_{j=1}^{k_r}|\hat p_{n,j}-p_{r,j}^*|$,
where $\hat p_n$ is the empirical cell-frequency vector.  Under
Model~C this quantity converges to zero at the predicted
$O\{(k_r/n)^{1/2}\}$ rate.  Under Model~M it stabilizes near $0.322$,
the fixed total-variation distance $\TV(q_{M,r},p_r^*)$ between the
Model-M cell law and the projected baseline $p_r^*=q_r$, the
pseudo-true projection of the parametric family at $\theta_0=0$.
The projection of the Model-M law itself is a distinct object, with
$\TV(q_{M,r},p_{M,r}^*)\approx 0.326$.

\begin{figure}[!t]
  \centering
  \includegraphics[width=0.52\textwidth]{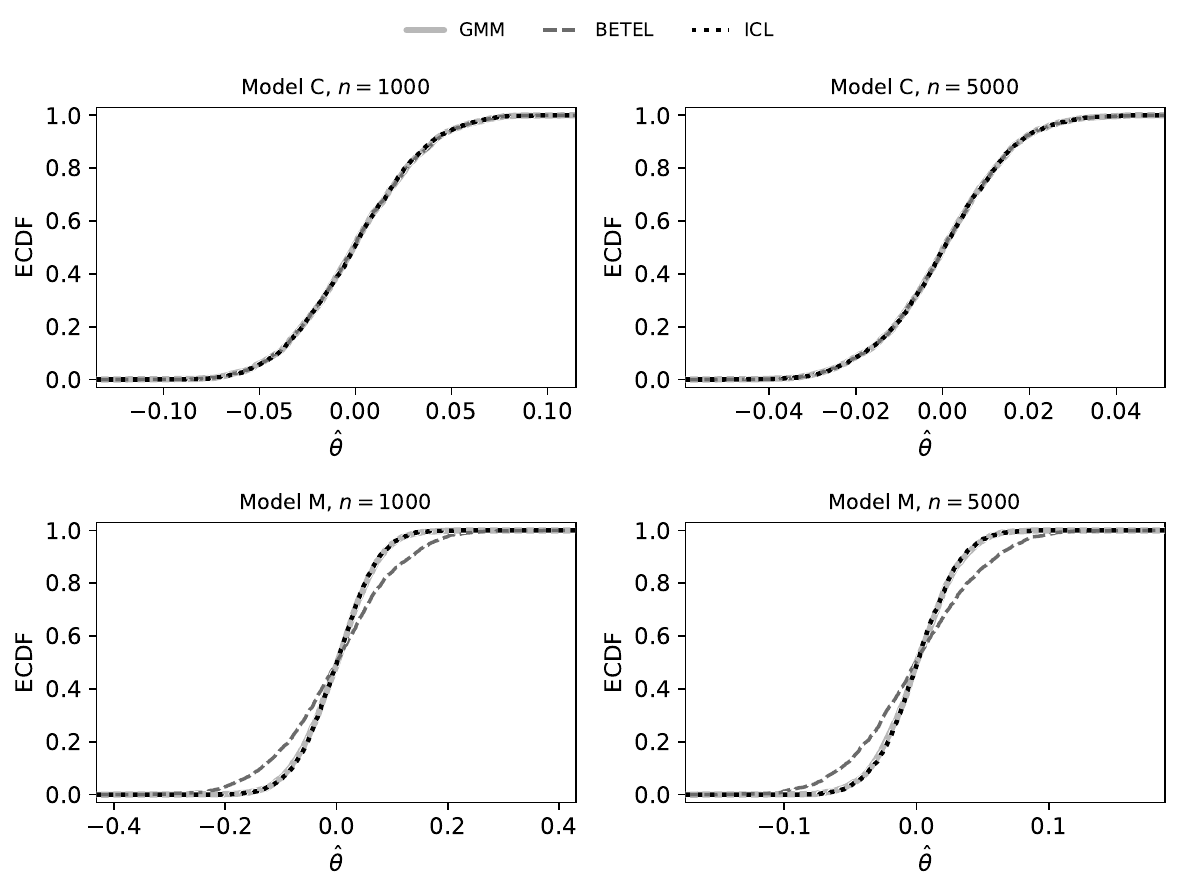}
  \figwidth=\textwidth
  \caption{Empirical distribution functions of $\hat\theta$ over $2000$
    replications.  Upper row: Model~C; lower row: Model~M ($\sigma=2$).
    Under correct specification the three procedures overlap; under
    misspecification $\hat\theta_{\text{\textsc{betel}}}$ is more dispersed
    while $\hat\theta_{\text{\textsc{icl}}}$ stays close to
    $\hat\theta_{\text{\textsc{gmm}}}$}
  \label{fig:sampling-supp}
  \alttext{In each panel the three empirical distribution functions overlap
  under correct specification.  Under misspecification the tilted-sample
  estimator is more dispersed while the other two nearly overlap.}
\end{figure}

\subsection{Outputs}

The parametric simulation script produces
five figures (geometry projection, curvature refinement, sampling
comparison, empirical cell discrepancy, standard-deviation comparison) and
two tables (main results, curvature under refinement).  The
\textsc{Python} script \texttt{simulation\_final.py} accompanies the
resubmission as a self-contained simulation package; the direct
conditional calculation is reproduced by
\texttt{direct\_block\_collapse.py}.

\section{Method Comparison}\label{sec:comparison}

The comparison with Bayesian asymptotics based on the exponentially
tilted empirical likelihood begins with a
difference in primitives.  In the present paper the
finite chart is derived from a known permutation-invariant constraint
via the empirical-measure reduction, whereas in the approach of
\citet{chibsimoni2018,chibsimoni2022} the tilted empirical likelihood on
the finite sample-weight polytope is taken as the inferential starting
point.
Theorem~A1 relates the ambient and discretized experiments at fixed sample
size.  Our emphasis is finite-sample geometry, rather than another
derivation of Bernstein--von~Mises theory.

Table~\ref{tab:comparison} summarizes the comparison.

\begin{table}[t]
\centering
\caption{Comparison of primitives and emphases across related
approaches to moment-conditioned inference.  GMM, generalized method of
moments; ETEL/BETEL, (Bayesian) exponentially tilted empirical
likelihood; BvM, Bernstein--von~Mises}\label{tab:comparison}
\footnotesize
\setlength{\tabcolsep}{4pt}
\begin{tabular}{p{1.8cm} p{2.7cm} p{2.7cm} p{2.5cm} p{2.7cm}}
\toprule
Method & Primitive object & How the finite polytope enters & Main
asymptotics & Finite-sample geometry \\
\midrule
GMM
& Quadratic discrepancy in sample moments
& Not central; the empirical moment vector is primary
& Classical asymptotic normality
& No explicit predictive-collapse or curvature analysis \\[4pt]

ETEL / BETEL
& Exponentially tilted empirical likelihood on sample weights
& The simplex of reweighted empirical probabilities is the inferential
carrier
& Likelihood-type asymptotics via tilted empirical likelihood
& Not focused on finite-sample localization \\[4pt]

Chib--Shin--Simoni
& Bayesian ETEL/BETEL-type likelihood
& The sample-weight polytope is taken as the primitive likelihood device
& BvM, misspecification, and model comparison under general conditions
& Finite polytope present, but geometric bounds are not the focus \\[4pt]

Present paper
& Fixed reference law $Q$ + known permutation-invariant constraint; optional
prior on $\theta$ at the inferential layer
& The finite polytope is \emph{derived} from the empirical measure via a
partition chart
& Le~Cam transfer at the experiment and projected-law level; chart LAN and
misspecification geometry, with broader BvM theory due to Chib--Shin--Simoni
& Yes: localization, reduced-Hessian curvature, and
predictive-collapse bounds \\
\bottomrule
\end{tabular}
\end{table}

\section{Partition-Sizing Experiment}\label{sec:partition-sizing}

To make the bias--variance trade-off in Remark~6 of the main paper visible,
we use the mean-and-variance moment functions of
\citet{schennach2007el},
$g(x,\theta)=\{x-\theta,(x-\theta)^2-1\}$ with baseline $Q=N(0,1)$
and a shifted target $\alpha=(0{.}5,\,1{.}25)$, so that the true
constrained law is $P^*=N(0{.}5,\,1)$ and the information projection
involves a genuine exponential tilt.

\begin{figure}[!t]
\centering
\begin{minipage}[t]{0.54\textwidth}
\centering
\includegraphics[width=\linewidth]{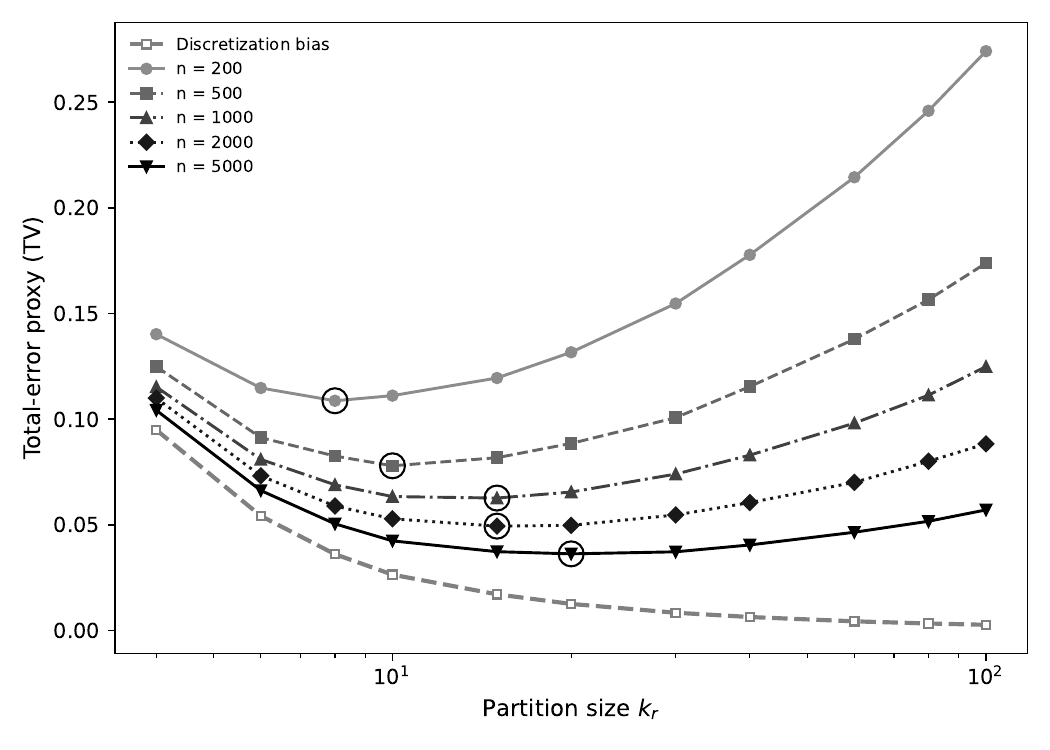}
\end{minipage}\hfill
\begin{minipage}[t]{0.44\textwidth}
\centering
\includegraphics[width=\linewidth]{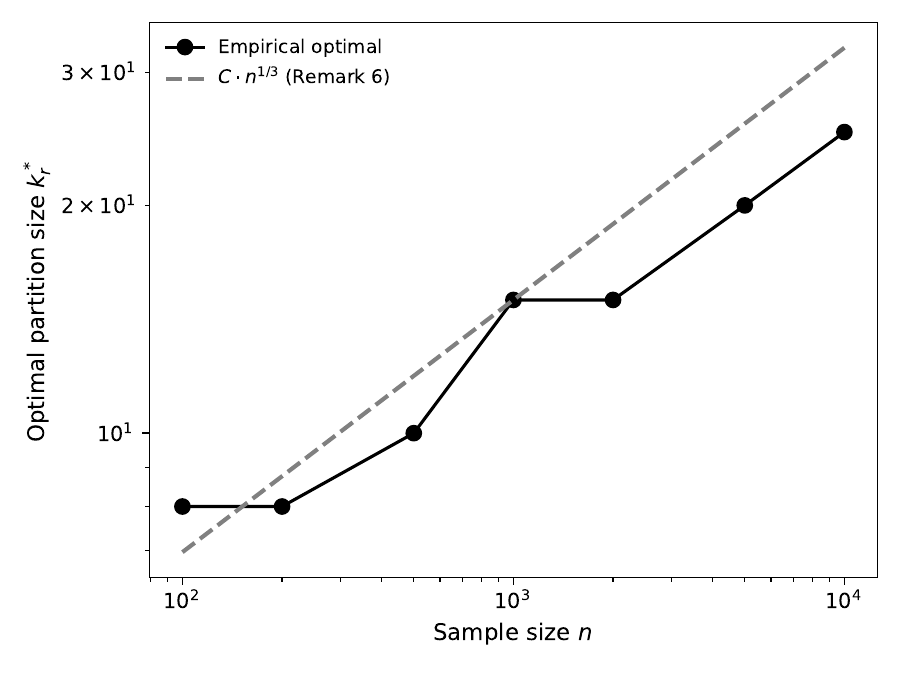}
\end{minipage}
\figwidth=\textwidth
\caption{Partition-sizing experiment.  Left: total-error proxy
  (chart-approximation bias + sampling fluctuation) versus partition size
  $k_r$, for $n\in\{200,500,1000,2000,5000\}$; the dashed grey curve is
  the discretization bias alone and open circles mark the minimizing
  $k_r$ for each $n$.  Right: optimal partition size $k_r^*$ versus
  sample size $n$ on a log-log scale; the dashed line is the $C\,n^{1/3}$
  reference from the heuristic of Remark~6; the empirical scaling
  exponent is $0{.}26$.}
\label{fig:partition-experiment}
\alttext{Left, U-shaped total-error curves over partition size have minima
that move right as the sample size grows.  Right, the selected optimal
partition sizes rise approximately linearly on a log-log scale.}
\end{figure}

\paragraph{Experiment~1: bias--variance trade-off.}
For each partition size $k_r\in\{4,6,\ldots,100\}$ and sample size
$n\in\{200,500,1000,2000,5000\}$, we compute two quantities:
the \emph{discretization bias}, defined as
$\TV(p_r^*,\,p_{\mathrm{true}})$ where $p_{\mathrm{true}}$ is the
vector of $P^*$-cell probabilities; and the \emph{sampling error},
defined as the mean $\TV(\widehat{p}_n,\,p_{\mathrm{true}})$ over
$300$ replications of data drawn from $P^*$.
The left panel of Fig.~\ref{fig:partition-experiment} plots their sum,
$\TV_{\mathrm{disc}}+\TV_{\mathrm{samp}}$, against $k_r$.  This
total-error proxy combines chart approximation and sampling fluctuation;
it is not the conditional predictive law.  Its U shape has the expected
interpretation: coarse partitions retain discretization bias, whereas very
fine ones pay more sampling error.  Open circles mark the minima, which
move to the right as $n$ grows and finer partitions become affordable.

\paragraph{Experiment~2: optimal partition scaling.}
For each $n\in\{100,200,\ldots,10000\}$ we select $k_r^*$ as the
minimizer of total error over
$k_r\in\{4,6,8,10,12,15,20,25,30,40,50,60,80,100\}$.
The right panel of Fig.~\ref{fig:partition-experiment} plots
$k_r^*$ against $n$ on a log-log
scale.  A linear fit gives $\log k_r^*=0{.}78+0{.}262\log n$ (empirical
exponent $0{.}26$), below the $1/3$ heuristic of Remark~6 for Lipschitz
moment maps on $\R^1$.  The candidate grid is discrete and the objective is
shallow near its minimum, so this exponent is descriptive.  The total-error
proxy at the selected partition decays as $n^{-0{.}336}$, close to
$O(n^{-1/3})$; neither fit is a proved uniform theorem.

\paragraph{Experiment~3: sensitivity to partition choice.}
In the correctly specified Schennach design ($P_0=Q=N(0,1)$,
$\alpha=(0,1)$), we vary $k_r\in\{6,8,\ldots,80\}$ and compute the
Monte Carlo standard deviation of the induced conditional-likelihood
estimator over $300$ replications at $n=1000$.  The standard deviation is
stable across partitions, ranging over $3{.}06$--$3{.}39\times 10^{-2}$
and staying close to the parametric rate $1/\sqrt{n}$.

The simulation code is available as
\texttt{simulation\_partition\_sizing.py} in the simulation package.

\end{document}